\documentclass[11pt,fleqn]{article}
\usepackage{stmaryrd}
\usepackage{amssymb}
\usepackage{bbm}
\usepackage{mathrsfs}
\usepackage{amsfonts}
\usepackage{amsmath,epsf,eufrak}
\begin{document}
\renewcommand{\theequation}{\arabic {section}.\arabic {equation}}
\newcommand{\bi}{\begin{equation}}
\newcommand{\ei}{\end{equation}}
\date{}
\title{  On a variance for twins of $k-$free numbers in arithmetic
progressions
}
\author{Zaizhao MENG}
\maketitle \baselineskip 24pt
 \begin{center}
 In this paper, we give a new upper bound of Barban-Davenport-Halberstam type for
  twins of  $k-$free numbers in arithmetic progressions.
\end{center}
\section{Introduction}
\setcounter{equation}{0}
     In recent years, Br\"{u}dern and others [1,2] have made a breakthrough in the circle method.
They give successful treatment, through the circle method, of binary
additive problems involving $k-$free numbers. Their results depend
upon the variance for $k-$free numbers (or twins of $k-$free
numbers) in arithmetic
progressions.\\
    In this paper, we give new results on the variance for twins of $k-$free numbers in arithmetic
progressions. Such result is analogous to the Barban-Davenport-Halberstam theorem for the primes in
arithmetic progressions.\\
    Let $\mu_{k}(n)$ be the characteristic function of the $k-$free numbers,
\bi \mu_{k}(n)=\sum\limits_{d^{k}\mid n}\mu(d), \ei
where $\mu(n)$ is the M\"{o}bius function.\\
Let \\
\bi A_{k}(x;q,a)= \sum\limits_{\stackrel{n\leq x}{n\equiv a(mod\
q)}}\mu_{k}(n)\mu_{k}(n+1)
 \ei
\bi g(q,a)= \sum^{\infty}\limits_{\stackrel{u,v=1}{(u^{k},q)\mid a,
(v^{k},q)\mid a+1}}\frac{\mu(uv)}{u^{k}v^{k}}(q,u^{k}v^{k})
 \ei
By Lemma 3.1 in [2], we have
$$ A_{k}(x;q,a)=q^{-1}g(q,a)x+O(x^{\frac{2}{k+1}+\varepsilon}).$$
We consider the variance \bi Y_{k}(x,Q)=\sum\limits_{q\leq
Q}\sum\limits_{a=1}^{q}\mid A_{k}(x;q,a)-q^{-1}g(q,a)x\mid^{2}.
 \ei
In [2], Br\"{u}dern,  Perelli and  Wooley obtained \bi
 Y_{k}(x,Q)\ll
 x^{\frac{2}{k}+\varepsilon}Q^{2-\frac{2}{k}}+x^{\frac{4}{k+1}+\varepsilon},\ \ \ \ 1<Q\leq x .
 \ei
In this paper, we obtain the following\\
{\bf THEOREM.}\  Suppose that $1< Q\leq x$ and that $k$ is an
integer with $k\geq 2$. Then
$$
 Y_{k}(x,Q)\ll
 x^{\frac{1}{k}+\varepsilon}Q^{2-\frac{1}{k}-\varepsilon}+
 x^{1+\frac{2}{k}}\log Q+x^{\frac{3}{2}+\frac{1}{2k}+\varepsilon}.
$$
Actually, the first term is
$Q^{2}(\frac{x}{Q})^{\frac{1}{k}+\varepsilon}$. For $k>2$, this
result is superior to (1.5) when $Q\gg x^{\frac{3}{4}}$.\\
We have
 \bi
 Y_{k}(x,Q)=S_{1}(x,Q)-2xS_{2}(x,Q)+x^{2}S_{3}(x,Q)
 \ei
where
 \bi
 S_{1}(x,Q)=\sum\limits_{q\leq Q}\sum\limits_{a=1}^{q}
A_{k}^{2}(x;q,a), \ei
\bi
 S_{2}(x,Q)=\sum\limits_{q\leq Q}\sum\limits_{a=1}^{q}
A_{k}(x;q,a)q^{-1}g(q,a),
 \ei
\bi
 S_{3}(x,Q)=\sum\limits_{q\leq Q}\sum\limits_{a=1}^{q}
q^{-2}g(q,a)^{2}.
 \ei
We notice that the function $g(q,a)$ does not depend only on $q$ and
$(a,q)$, unlike most sequences investigated before.\\
In section 2-3, we discuss $S_{2}(x,Q)$ and $S_{3}(x,Q)$, these
depend on the solutions of the corresponding congruences. In section
4-7, we discuss $S_{1}(x,Q)$ by the Hardy-Littlewood method on the
line of Vaughan[9]. In section 8, we discuss the singular series
$\mathfrak{S}(n)$. In section 9-10, we use the Hurwitz Zeta function
to discuss the generating function of $\mathfrak{S}(uv)$. In section
11, we use a trick to avoid the difficulty of calculating some
constants and then we obtain the Theorem.\\
 {\bf Notation:} Let $k>1$ denote a positive integer. Throughout,
$\varepsilon$ is a sufficiently small positive number, the implicit
constants in Vinogradov's notation $\ll$, and in Landau's
O-notation, will depend at most on $k,\varepsilon$ unless it is
pointed out depend upon the corresponding parameters. $n\equiv
a(mod\ q)$ may be written as $n\equiv a(q)$.\ The greatest common
divisor and the least common multiple of integers $a,b$ are denoted
by $(a,b)$ and $[a,b]$ respectively; $\mu(n)$ denotes the M\"{o}bius
function and $\tau(n)$ denotes the divisor function; $[w]$ denotes
the integer part of $w$ and $e(\alpha)=\exp(2\pi i \alpha)$;
$\sum\limits_{a=1}^{q}{}^{^{\prime}}$ means $(a,q)=1$. The letter
$p$ denotes a prime number, and write $p^{t}\parallel n$ when
$p^{t}\mid n$ but $p^{t+1}\dagger n$.  Let $x$ denote a sufficiently
large real number and $Q$ be a positive real number with $1<Q\leq
x$.
\section{ The formula for $S_{2}(x,Q)$ }
\setcounter{equation}{0}
For fixed positive integers $u,v,q,r,s$,
let \bi
 J_{1}=
\sum^{q}\limits_{\stackrel{a=1}{(u^{k},q)\mid a, (v^{k},q)\mid
a+1}}\sum\limits_{\stackrel{n\leq x,n\equiv a(mod\ q)}{r^{k}\mid n,
s^{k}\mid n+1}}1,\ \ \
s_{0}=(\frac{(u^{k},q)}{(r^{k},u^{k},q)},\frac{r^{k}}{(r^{k},q)}).
\ei
{\bf LEMMA 2.1.}We have
$$
J_{1}=\frac{(r^{k},u^{k},q)s_{0}}{r^{k}s^{k}(u^{k},q)(v^{k},q)}
(s^{k},\frac{r^{k}(u^{k},q)(v^{k},q)}{(r^{k},u^{k},q)})x+O(1),
$$
provided that
\bi
(\frac{(u^{k},q)(r^{k},q)}{(r^{k},u^{k},q)},(v^{k},q))=(\frac{(u^{k},q)}{(r^{k},u^{k},q)s_{0}},s^{k})=(r,s)=(q,u,v)=1.
\ei
Proof.\ We use an idea of [4]. Let $n=r^{k}b,\ a=(u^{k},q)c$, then\\
\bi
 J_{1}= \sum^{q(u^{k},q)^{-1}}\limits_{\stackrel{c=1}{(v^{k},q)\mid
(u^{k},q)c+1}}\sum\limits_{\stackrel{r^{k}b\leq x}{r^{k}b\equiv
  (u^{k},q)c(mod\ q),s^{k}\mid r^{k}b+1}}1,\ \ \ \ (r,s)=(q,u,v)=1.\ \ \
\ei
Now,\\
$\frac{r^{k}}{(r^{k},q)}b\equiv
  \frac{(u^{k},q)c}{(r^{k},q)}(mod\ \frac{q}{(r^{k},q)}), \ \ (r^{k},q)\mid
  (u^{k},q)c.$\\
  Since  $
\frac{(u^{k},q)(r^{k},q)}{(r^{k},u^{k},q)}=[(u^{k},q),(r^{k},q)]\mid
q ,$  we write $c=\frac{(r^{k},q)}{(r^{k},u^{k},q)}d$, then\\
$
 J_{1}=\#\{1\leq b\leq xr^{-k},\ 1\leq d\leq
\frac{q(r^{k},u^{k},q)}{(u^{k},q)(r^{k},q)}:\ (v^{k},q)\mid
 (u^{k},q)\frac{(r^{k},q)}{(r^{k},u^{k},q)}d+1,\\
 \frac{r^{k}}{(r^{k},q)}b\equiv
  \frac{(u^{k},q)}{(r^{k},q)}\frac{(r^{k},q)}{(r^{k},u^{k},q)}d(mod\
    \frac{q}{(r^{k},q)}),\ s^{k}\mid r^{k}b+1.\}
$\\
We have\\
$$
\frac{r^{k}}{(r^{k},q)}b\equiv
  \frac{(u^{k},q)d}{(r^{k},u^{k},q)}(mod\
    \frac{q}{(r^{k},q)}),
$$
hence \\
$
\frac{r^{k}}{(r^{k},q)}\frac{(r^{k},u^{k},q)}{(u^{k},q)}b\equiv
  d(mod\
    \frac{q}{(r^{k},q)}\frac{(r^{k},u^{k},q)}{(u^{k},q)}),
$\\
and \bi
\frac{(u^{k},q)}{(r^{k},u^{k},q)}\mid\frac{r^{k}}{(r^{k},q)}b.
\ei
Write $b=\frac{(u^{k},q)}{(r^{k},u^{k},q)s_{0}}b_{1}$,\ then\\
$
 J_{1}=\#\{1\leq \frac{(u^{k},q)}{(r^{k},u^{k},q)s_{0}}b_{1}\leq xr^{-k},\ 1\leq d\leq
\frac{q(r^{k},u^{k},q)}{(u^{k},q)(r^{k},q)}:\ (v^{k},q)\mid
 (u^{k},q)\frac{(r^{k},q)}{(r^{k},u^{k},q)}d+1,\\
 \frac{r^{k}}{(r^{k},q)s_{0}}b_{1}\equiv d
  (mod\
    \frac{q(r^{k},u^{k},q)}{(r^{k},q)(u^{k},q)}),\ s^{k}\mid r^{k}\frac{(u^{k},q)}{(r^{k},u^{k},q)s_{0}}b_{1}+1.\}
$\\
We have \\
\bi
 \frac{(u^{k},q)(r^{k},q)}{(r^{k},u^{k},q)}d\equiv -1(mod\
 (v^{k},q)),
 \ei
 and
 \bi
 (\frac{(u^{k},q)(r^{k},q)}{(r^{k},u^{k},q)},(v^{k},q))=1,
 \ei
 \bi
 (v^{k},q)\mid\frac{q(r^{k},u^{k},q)}{(u^{k},q)(r^{k},q)},
 \ei
\bi (\frac{(u^{k},q)}{(r^{k},u^{k},q)s_{0}},s^{k})=1.
 \ei
 We write \\
 $
d=d_{1}+m(v^{k},q),\ 1\leq d_{1}\leq (v^{k},q),\ 0\leq
m<\frac{q(r^{k},u^{k},q)}{(v^{k},q)(u^{k},q)(r^{k},q)},
 $\\
 then
 \bi
\frac{r^{k}}{(r^{k},q)s_{0}}b_{1}\equiv d_{1}+m(v^{k},q)
  (mod\
    \frac{q(r^{k},u^{k},q)}{(r^{k},q)(u^{k},q)})
 \ei
\bi
 \frac{r^{k}}{(r^{k},q)s_{0}}b_{1}\equiv d_{1} (mod\ (v^{k},q)).
 \ei
If a prime $p\mid (v,q)$, then $p\mid v,\ p\mid q, p\dagger u$,
again if $p\mid\frac{r^{k}}{(r^{k},q)s_{0}}$, then $p\mid r,$
hence,by (2.6)
$p\mid(\frac{(u^{k},q)(r^{k},q)}{(r^{k},u^{k},q)},(v^{k},q))=1,$
  this is a contradiction.\\
   We deduce that
 \bi
((v^{k},q),\frac{r^{k}}{(r^{k},q)s_{0}})=1.
 \ei
Now the congruence (2.10) has a unique solution modulus $(v^{k},q)$.\\
Let $b_{1}=b_{0}+t(v^{k},q), \ 1\leq b_{0}\leq (v^{k},q)$, by (2.7)
and (2.9)
 \bi \frac{r^{k}(v^{k},q)}{(r^{k},q)s_{0}}t\equiv
d_{1}-\frac{r^{k}}{(r^{k},q)s_{0}}b_{0}+m(v^{k},q)
  (mod\ \frac{q(r^{k},u^{k},q)}{(r^{k},q)(u^{k},q)}).
\ei
Also, we have
\bi
 \frac{r^{k}(u^{k},q)}{(r^{k},u^{k},q)s_{0}}(v^{k},q)t\equiv
 -1-\frac{r^{k}(u^{k},q)}{(r^{k},u^{k},q)s_{0}}b_{0} \ (mod\ s^{k}).
\ei
By (2.7), (2.10) and (2.12),
\bi
\frac{r^{k}}{(r^{k},q)s_{0}}t\equiv
(d_{1}-\frac{r^{k}}{(r^{k},q)s_{0}}b_{0})(v^{k},q)^{-1}+m
  (mod\ \frac{q(r^{k},u^{k},q)}{(r^{k},q)(u^{k},q)(v^{k},q)}),
\ei
 and
\bi
 d_{1}\equiv \frac{r^{k}}{(r^{k},q)s_{0}}b_{0}\ (mod\ (v^{k},q)).
\ei
By (2.5),
\bi
 \frac{(u^{k},q)(r^{k},q)}{(r^{k},u^{k},q)}d_{1}\equiv -1(mod\
 (v^{k},q)).
 \ei
By (2.14), we only want to count the number of $t$ satisfying (2.13)
and
$$
\frac{(u^{k},q)}{(r^{k},u^{k},q)s_{0}}(b_{0}+t(v^{k},q))\leq
xr^{-k},
$$
that is
\bi
t\leq\frac{x(r^{k},u^{k},q)s_{0}}{r^{k}(u^{k},q)(v^{k},q)}-\frac{b_{0}}{(v^{k},q)}.
 \ei
By (2.5), (2.15) and (2.16),
\bi
 \frac{(u^{k},q)r^{k}}{(r^{k},u^{k},q)s_{0}}b_{0}\equiv -1(mod\
 (v^{k},q)).
 \ei
 Consequently, by(2.13) and (2.18),
 \bi
 ((v^{k},q),s^{k})\mid 1+\frac{(u^{k},q)r^{k}}{(r^{k},u^{k},q)s_{0}}b_{0}.
  \ei
  By (2.3), (2.8) and (2.19), we have
  $$
(s^{k},\frac{r^{k}(u^{k},q)(v^{k},q)}{(r^{k},u^{k},q)s_{0}})\mid
1+\frac{(u^{k},q)r^{k}}{(r^{k},u^{k},q)s_{0}}b_{0}.
  $$
By (2.13) and $(s,s_{0})=1$, we obtain\\
$
 \frac{r^{k}(u^{k},q)(v^{k},q)}{(r^{k},u^{k},q)}(s^{k},\frac{r^{k}(u^{k},q)(v^{k},q)}{(r^{k},u^{k},q)})^{-1}t\\
 \equiv
 (-s_{0}-\frac{r^{k}(u^{k},q)}{(r^{k},u^{k},q)}b_{0})(s^{k},\frac{r^{k}(u^{k},q)(v^{k},q)}{(r^{k},u^{k},q)})^{-1}
 \ (mod\ s^{k}(s^{k},\frac{r^{k}(u^{k},q)(v^{k},q)}{(r^{k},u^{k},q)})^{-1}).
$\\
Therefore, the number of $t$ is
$$
\frac{x(r^{k},u^{k},q)s_{0}}{r^{k}(u^{k},q)(v^{k},q)}(s^{k}(s^{k},\frac{r^{k}(u^{k},q)(v^{k},q)}{(r^{k},u^{k},q)})^{-1})^{-1}
+O(1),
$$
that is
$$
J_{1}=\frac{(r^{k},u^{k},q)s_{0}}{r^{k}s^{k}(u^{k},q)(v^{k},q)}
(s^{k},\frac{r^{k}(u^{k},q)(v^{k},q)}{(r^{k},u^{k},q)})x+O(1).
$$
The lemma follows.\\
Now,we give the formula for $S_{2}(x,Q)$.\\
$ S_{2}(x,Q)=\sum\limits_{q\leq
Q}q^{-1}\sum^{\infty}\limits_{u,v=1}\frac{\mu(uv)}{u^{k}v^{k}}(q,u^{k}v^{k})\sum\limits_{\stackrel{a=1}{(u^{k},q)\mid
a, (v^{k},q)\mid a+1}}^{q} \sum\limits_{\stackrel{n\leq x}{n\equiv
a(mod\ q)}}\mu_{k}(n)\mu_{k}(n+1)\\
= \sum\limits_{q\leq
Q}q^{-1}\sum^{\infty}\limits_{u,v=1}\frac{\mu(uv)}{u^{k}v^{k}}(q,u^{k}v^{k})\sum\limits_{\stackrel{r^{k}\leq
x}{s^{k}\leq x+1}}\mu(r)\mu(s)J_{1}.
 $\\
By Lemma 2.1,
\bi
S_{2}(x,Q)=S_{2A}(x,Q)+S_{2B}(x,Q),
 \ei
where\\
$ S_{2A}(x,Q)=\sum\limits_{q\leq
Q}q^{-1}\sum^{\infty}\limits_{u,v=1}\frac{\mu(uv)}{u^{k}v^{k}}(q,u^{k}v^{k})
\sum\limits_{\stackrel{r^{k}\leq x,s^{k}\leq
x+1}{(2.1),(2.2)}}\mu(r)\mu(s)\frac{(r^{k},u^{k},q)s_{0}}{r^{k}s^{k}(u^{k},q)(v^{k},q)}
(s^{k},\frac{r^{k}(u^{k},q)(v^{k},q)}{(r^{k},u^{k},q)})x, $\\
$ S_{2B}(x,Q)\ll\sum\limits_{q\leq
Q}q^{-1}\sum^{\infty}\limits_{u,v=1}\frac{\mu(uv)}{u^{k}v^{k}}(q,u^{k}v^{k})\sum\limits_{\stackrel{r^{k}\leq
x}{s^{k}\leq x+1}}1, $\\
and $\sum\limits_{\stackrel{r^{k}\leq x,s^{k}\leq
x+1}{(2.1),(2.2)}}$ means $u,v,r,s,q$ satisfy (2.1) and (2.2).\\
We give upper bound for $S_{2B}(x,Q)$ first.\\
$ S_{2B}(x,Q) \ll
x^{\frac{2}{k}}\sum^{\infty}\limits_{m=1}\frac{\mu^{2}(m)\tau(m)}{m^{k}}\sum\limits_{q\leq
Q}q^{-1}(q,m^{k})\\
\ll
x^{\frac{2}{k}}\sum^{\infty}\limits_{m=1}\frac{\mu^{2}(m)\tau(m)}{m^{k}}\sum\limits_{t\mid
m^{k}}\sum\limits_{\stackrel{q\leq Q}{(q,m^{k})=t}}q^{-1}(q,m^{k})
\ll
x^{\frac{2}{k}}\sum^{\infty}\limits_{m=1}\frac{\mu^{2}(m)\tau(m)}{m^{k}}\sum\limits_{t\mid
m^{k}}\sum\limits_{qt\leq
Q}q^{-1}$,\\
hence
 \bi
 S_{2B}(x,Q)\ll x^{\frac{2}{k}}\log Q.
 \ei
By
$(\frac{(u^{k},q)}{(r^{k},u^{k},q)},\frac{r^{k}}{(r^{k},q)}\frac{(r^{k},q)}{(r^{k},u^{k},q)})
=(\frac{(u^{k},q)}{(r^{k},u^{k},q)},\frac{r^{k}}{(r^{k},u^{k},q)})=1,$
  we have
 \bi
 s_{0}=(\frac{(u^{k},q)}{(r^{k},u^{k},q)},\frac{r^{k}}{(r^{k},q)})=1.
 \ei
 Since $(u,v)=1$, the conditions (2.2) and (2.22) are equivalent to
\bi (r,v,q)=(s,u,q)=(s,r)=1.
 \ei
By (2.2),
 \bi
  S_{2A}(x,Q)=x\sum\limits_{q\leq
Q}q^{-1}\sum^{\infty}\limits_{u,v=1}\frac{\mu(uv)}{u^{k}v^{k}}
\sum\limits_{\stackrel{r^{k}\leq x,s^{k}\leq
x+1}{(2.23)}}\mu(r)\mu(s)r^{-k}s^{-k}(s^{k},v^{k},q)(r^{k},u^{k},q).
 \ei
The summations of $r,s$ can be completed to $\infty$.\\
\bi
  S_{2A}(x,Q)=c_{1}(Q)x+S_{2A_{0}},
\ei
where
\bi c_{1}(Q)= \sum\limits_{q\leq
Q}q^{-1}\sum^{\infty}\limits_{u,v=1}\frac{\mu(uv)}{u^{k}v^{k}}
\sum^{\infty}\limits_{\stackrel{r,s=1}{(2.23)}}\mu(r)\mu(s)r^{-k}s^{-k}(s^{k},v^{k},q)(r^{k},u^{k},q),
 \ei
and\\
$ S_{2A_{0}}\ll x\sum\limits_{q\leq
Q}q^{-1}\sum^{\infty}\limits_{u,v=1}\frac{1}{u^{k}v^{k}}
(\sum\limits_{\stackrel{r^{k}>x}{s^{k}\leq
x+1,(r,s)=1}}r^{-k}s^{-k}(s^{k},v^{k},q)(r^{k},u^{k},q) +\\
\sum\limits_{\stackrel{r\geq 1}{s^{k}>
x+1,(r,s)=1}}r^{-k}s^{-k}(s^{k},v^{k},q)(r^{k},u^{k},q))
=S_{2Aa}+S_{2Ab},\ \
 $say.\\
$ S_{2Aa} \ll x\sum\limits_{q\leq
Q}q^{-1}\sum^{\infty}\limits_{u,v=1}\frac{1}{u^{k}v^{k}}
\sum\limits_{s^{k}\leq x+1}s^{-k}(s^{k},v^{k},q)\sum\limits_{t\mid
u}\sum\limits_{\stackrel{r^{k}>x}{(r,u)=t,(r,s)=1}}r^{-k}(t^{k},q)\\
 \ll x\sum\limits_{q\leq
Q}q^{-1}\sum^{\infty}\limits_{u,v=1}\frac{1}{u^{k}v^{k}}
\sum\limits_{s^{k}\leq x+1}s^{-k}(s^{k},v^{k},q)\sum\limits_{t\mid
u,(t,s)=1}(t^{k},q)t^{-k}\sum\limits_{t^{k}r^{k}>x}r^{-k}.$\\
Now\\
 $
\sum\limits_{t\mid
u,(t,s)=1}(t^{k},q)t^{-k}\sum\limits_{t^{k}r^{k}>x}r^{-k}\\
 \ll
\sum\limits_{\stackrel{t\mid u,(t,s)=1}{t^{k}\leq
x}}(t^{k},q)t^{-k}(x^{\frac{1}{k}}t^{-1})^{1-k} +
\sum\limits_{\stackrel{t\mid u}{t^{k}> x}}(t^{k},q)t^{-k}
 \ll  x^{\frac{1}{k}-1}\sum\limits_{\stackrel{t\mid u,(t,s)=1}{t^{k}\leq
x}}(t^{k},q)t^{-1} + \sum\limits_{\stackrel{t\mid u,(t,s)=1}{t^{k}>
x}}(t^{k},q)t^{-k}.
 $\\
Hence\\
$ S_{2Aa} \ll x\sum\limits_{q\leq
Q}q^{-1}\sum^{\infty}\limits_{u,v=1}\frac{1}{u^{k}v^{k}}
\sum\limits_{s^{k}\leq
x+1}s^{-k}(s^{k},v^{k},q)(x^{\frac{1}{k}-1}\sum\limits_{\stackrel{t\mid
u,(t,s)=1}{t^{k}\leq x}}(t^{k},q)t^{-1} +
\sum\limits_{\stackrel{t\mid u,(t,s)=1}{t^{k}>
x}}(t^{k},q)t^{-k})\\
\ll S_{2Aa1}+S_{2Aa2},$  \ say.\\
Firstly, we estimate  $ S_{2Aa1}$,\\
$ S_{2Aa1}\ll
x^{\frac{1}{k}}\sum^{\infty}\limits_{u,v=1}\frac{1}{u^{k}v^{k}}
\sum\limits_{s^{k}\leq x+1}s^{-k}\sum\limits_{\stackrel{t\mid
u,(t,s)=1}{t^{k}\leq x}}t^{-1}\sum\limits_{q\leq
Q}q^{-1}(s^{k},v^{k},q)(t^{k},q)\\
 \ll
x^{\frac{1}{k}}\sum^{\infty}\limits_{u,v=1}\frac{1}{u^{k}v^{k}}
\sum\limits_{s^{k}\leq x+1}s^{-k}\sum\limits_{\stackrel{t\mid
u}{t^{k}\leq x}}t^{-1}\sum\limits_{q\leq Q}q^{-1}(s^{k}t^{k},q)\\
 \ll
x^{\frac{1}{k}}\sum^{\infty}\limits_{u,v=1}\frac{1}{u^{k}v^{k}}
\sum\limits_{s^{k}\leq x+1}s^{-k}\sum\limits_{\stackrel{t\mid
u}{t^{k}\leq x}}t^{-1}\sum\limits_{j\mid
s^{k}t^{k}}\sum\limits_{qj\leq Q}q^{-1}
\ll x^{\frac{1}{k}}\log Q,$\\
hence
\bi
 S_{2Aa1}\ll  x^{\frac{1}{k}}\log Q.
\ei
Secondly, as the estimate of $ S_{2Aa1}$, we have\\
$ S_{2Aa2} \ll x\sum\limits_{q\leq
Q}q^{-1}\sum^{\infty}\limits_{u,v=1}\frac{1}{u^{k}v^{k}}
\sum\limits_{s^{k}\leq x+1}s^{-k}(s^{k},v^{k},q)
\sum\limits_{\stackrel{t\mid u,(t,s)=1}{t^{k}>
x}}(t^{k},q)t^{-k}\\
\ll
x\sum^{\infty}\limits_{u,v=1}\frac{1}{u^{k}v^{k}}
\sum\limits_{s^{k}\leq x+1}s^{-k} \sum\limits_{\stackrel{t\mid
u,(t,s)=1}{t^{k}> x}}\sum\limits_{q\leq
Q}q^{-1}(s^{k},v^{k},q)(t^{k},q)t^{-k} \\
\ll x\log Q\sum^{\infty}\limits_{v=1}\frac{1}{v^{k}}
\sum\limits_{s^{k}\leq
x+1}s^{-k}\tau(s^{k})(x^{\frac{1}{k}})^{1-2k+\varepsilon},$\\
hence
\bi
 S_{2Aa2}\ll  x^{\frac{1}{k}-1+\varepsilon}\log Q.
\ei
By (2.27) and (2.28),
\bi
 S_{2Aa}\ll x^{\frac{1}{k}}\log Q+x^{\frac{1}{k}-1+\varepsilon}\log Q.
\ei In the same manner we have
\bi
 S_{2Ab}\ll x^{\frac{1}{k}}\log Q+x^{\frac{1}{k}-1+\varepsilon}\log Q.
\ei
By (2.29) and (2.30),
 \bi
 S_{2A_{0}}\ll x^{\frac{1}{k}}\log Q+x^{\frac{1}{k}-1+\varepsilon}\log Q.
\ei
By (2.20), (2.21), (2.25) and (2.31), we obtain \bi
 S_{2}(x,Q)=c_{1}(Q)x+O(x^{\frac{2}{k}}\log Q),
\ei
where $c_{1}(Q)$ is given by (2.26).
\section{ The formula for $S_{3}(x,Q)$ }
\setcounter{equation}{0}
By (1.3) and (1.9), \bi
 S_{3}(x,Q)=\sum\limits_{q\leq Q}\sum\limits_{a=1}^{q}
q^{-2} \sum^{\infty}\limits_{\stackrel{u,v=1}{(u^{k},q)\mid a,
(v^{k},q)\mid a+1}}\frac{\mu(uv)}{u^{k}v^{k}}(q,u^{k}v^{k})
\sum^{\infty}\limits_{\stackrel{r,s=1}{(r^{k},q)\mid a,
(s^{k},q)\mid a+1}}\frac{\mu(rs)}{r^{k}s^{k}}(q,r^{k}s^{k}).
 \ei
Let \bi
 \Delta =\#\{1\leq a\leq q:\ (u^{k},q)\mid a,\ (r^{k},q)\mid
a,\ (v^{k},q)\mid a+1 ,\   (s^{k},q)\mid a+1\}.
 \ei
Then
$$
(u,s,q)=(v,r,q)=1,
$$
and
$$
\Delta =\#\{1\leq m\leq \frac{q}{[(u^{k},q),(r^{k},q)]}:\
m[(u^{k},q),(r^{k},q)]+1\equiv 0(mod\ [(v^{k},q),(s^{k},q)])\},
$$
therefore,
$$
([(u^{k},q),(r^{k},q)],[(v^{k},q),(s^{k},q)])=1,
$$
and
\bi \Delta =\frac{q}{[(u^{k},q),(r^{k},q)][(v^{k},q),(s^{k},q)]}
=\frac{q(r^{k},u^{k},q)(s^{k},v^{k},q)}{(u^{k},q)(r^{k},q)(v^{k},q)(s^{k},q)}.
\ei Hence, by (3.1), (3.2) and (3.3),
 \bi
 S_{3}(x,Q)=c_{1}(Q),
\ei
where $c_{1}(Q)$ is given by (2.26).
\section{ Lemmas for $S_{1}(x,Q)$ }
\setcounter{equation}{0} As Vaughan had done in [9], our proof of
the theorem uses the Hardy-Littlewood method(see [8]) and depends
heavily on [2] where the bounds for
\bi
S(\alpha)=S_{k}(\alpha)=\sum\limits_{n\leq
x}\mu_{k}(n)\mu_{k}(n+1)e(n\alpha),
 \ei
and related expressions, are obtained.\\
Let $R=x^{1/2+\tau}$, where $\tau=\tau_{k}$ is a sufficiently small
positive number$(2\varepsilon<\tau<\frac{1}{10k})$, and let
$\mathfrak{M}$ denote the union of the intervals
 \bi
\mathfrak{M}(q,a)=\{\alpha :\mid q\alpha -a\mid\leq R^{-1}\},
\ei
with $1\leq a\leq q\leq x/R$ and $(q,a)=1$, and
$\mathscr{M}=(R^{-1},1+R^{-1}]\backslash \mathfrak{M}$.\\
 {\bf LEMMA 4.1.}We have
$$
\int_{\mathscr{M}}\mid S(\alpha)\mid^{2}d\alpha\ll
x^{\frac{1}{k}+\varepsilon}R^{1-\frac{1}{k}}+x^{2+\varepsilon}R^{\frac{2}{k}-3}+x^{\frac{4}{k+1}+1+\varepsilon}R^{-2}.$$
Proof.\ This is Theorem 2 of [2] by taking $Q=x/R$, we note that
$\mathfrak{M}\subseteq (R^{-1},1+R^{-1}]$, the interval
$(Q^{-1},1+Q^{-1}]$ of [2] may be replaced by $(R^{-1},1+R^{-1}]$.\\
Let \bi K(\alpha)=\sum\limits_{u\leq Q}\sum\limits_{v\leq
x/u}e(uv\alpha). \ei When $Q\leq \sqrt{x}$ define \bi
K_{q}(\alpha)=\sum\limits_{u\leq Q/q}\sum\limits_{v\leq
x/qu}e(quv\alpha), \ei
\bi H_{q}(\alpha)=\left\{\begin{array}{cc}
\sum\limits_{\stackrel{u\leq
Q}{q\dagger u}}\sum\limits_{v\leq x/u}e(uv\alpha),&(q>1),\\
0,&(q=1).
\end{array}
\right.
 \ei
When $Q>\sqrt{x}$ define \bi K_{q}(\alpha)=\sum\limits_{u\leq
\sqrt{x}/q}(\sum\limits_{v\leq
x/qu}e(quv\alpha)+\sum\limits_{\sqrt{x}<v\leq
\min(Q,x/qu)}e(quv\alpha)), \ei
\bi
H_{q}(\alpha)=\left\{\begin{array}{cc} \sum\limits_{\stackrel{u\leq
\sqrt{x}}{q\dagger u}}(\sum\limits_{v\leq x/u}e(uv\alpha)+\sum\limits_{\sqrt{x}<v\leq \min(Q,x/u)}e(uv\alpha)),&(q>1),\\
0,&(q=1).
\end{array}
\right.
 \ei
We have the following results on these functions(see the section 2
of Vaughan [9]). \\
{\bf LEMMA 4.2.} Suppose that $Q\leq x$ and $q\in N$. Then
$$K(\alpha)=K_{q}(\alpha)+ H_{q}(\alpha).$$
Proof.\ This is Lemma 2.9 of [9].\\
{\bf LEMMA 4.3.} Suppose that $(a,q)=1$ and $\mid q\alpha-a\mid\leq
q^{-1}$. Then
$$K(\alpha)\ll (xq^{-1}+q)\log x.$$
Proof.\ This is Lemma 2.10 of [9].\\
{\bf LEMMA 4.4.} Suppose that $(a,q)=1,q\in N$ and
$\alpha=a/q+\beta$. Then
$$K_{q}(\alpha)=K_{q}(\beta)\ll \sum\limits_{u\leq \sqrt{x}/q}\frac{x}{qu+x\|qu\beta\|}.$$
Moreover, if $\mid\beta\mid\leq 1/(2\sqrt{x})$, then
$$K_{q}(\alpha)\ll \frac{x\log x}{q+qx\mid\beta\mid}.$$
Proof.\ This is Lemma 2.11 of [9].\\
{\bf LEMMA 4.5.} Suppose that $(a,q)=1,q\in N$ and $\mid
q\alpha-a\mid\leq 1/(2\sqrt{x})$. Then
$$H_{q}(\alpha)\ll(\min(Q,\sqrt{x})+q)\log x.$$
Proof.\ This is Lemma 2.12 of [9].\\
By Dirichlet's theorem on diophantine approximation, Lemma 4.1 and
Lemma 4.3,
 \bi
 \sup\limits_{\mathscr{M}}\mid K(\alpha)\mid\ll R\log x
 \ei
and\\
$ \int_{\mathscr{M}}\mid S(\alpha)\mid^{2}K(\alpha)d\alpha\ll
x^{\frac{1}{k}+\varepsilon}R^{2-\frac{1}{k}}+x^{2+\varepsilon}R^{\frac{2}{k}-2}+x^{\frac{4}{k+1}+1+\varepsilon}R^{-1}
$\\
\bi \ll
x^{1+\frac{1}{k}-\frac{\tau}{2}}+x^{\frac{1}{2}+\frac{4}{k+1}
-\frac{\tau}{2}} . \ei
 By (1.2) and (1.7), \\
$
 S_{1}(x,Q)=\sum\limits_{q\leq Q}\sum\limits_{\stackrel{m,n\leq x}{m\equiv n(mod\
q)}}\mu_{k}(n)\mu_{k}(n+1)\mu_{k}(m)\mu_{k}(m+1)\\
=2\sum\limits_{q\leq Q}\sum\limits_{\stackrel{m<n\leq x}{m\equiv
n(mod\ q)}}\mu_{k}(n)\mu_{k}(n+1)\mu_{k}(m)\mu_{k}(m+1)
+[Q]\sum\limits_{n\leq x}\mu_{k}(n)\mu_{k}(n+1).$\\
Hence
\bi
 S_{1}(x,Q)=2S_{A}(x,Q)+[Q]\sum\limits_{n\leq
x}\mu_{k}(n)\mu_{k}(n+1) \ei
 where
 \bi
S_{A}(x,Q)=\sum\limits_{q\leq Q}\sum\limits_{\stackrel{m<n\leq
x}{m\equiv n(mod\ q)}}\mu_{k}(n)\mu_{k}(n+1)\mu_{k}(m)\mu_{k}(m+1).
 \ei
By (4.1), (4.3) and (4.11),
 \bi
 S_{A}(x,Q) =\int_{R^{-1}}^{1+R^{-1}}\mid S(\alpha)\mid^{2}K(\alpha)d\alpha
=\int_{\mathfrak{M}}\mid S(\alpha)\mid^{2}K(\alpha)d\alpha
+\int_{\mathscr{M}}\mid S(\alpha)\mid^{2}K(\alpha)d\alpha . \ei
 By
(4.9) and (4.12),
 \bi S_{A}(x,Q)=\int_{\mathfrak{M}}\mid
S(\alpha)\mid^{2}K(\alpha)d\alpha
+O(x^{1+\frac{1}{k}-\frac{\tau}{2}}+x^{\frac{1}{2}+\frac{4}{k+1}
-\frac{\tau}{2}}). \ei
By (1.3)
\bi
g(q,a)=\sum\limits_{n=1}^{\infty}\mu(n)\frac{(n^{k},q)}{n^{k}}\psi_{k}(n;q,a),
\ei
 where $\psi_{k}(n;q,a)$ denotes the number of pairs $u,v$ of
natural numbers with $uv=n,\ (u,v)=1,\ (u^{k},q)\mid a,\
(v^{k},q)\mid a+1$.\  For fixed $a,\ q$,\ the function
$\psi_{k}(n;q,a)$ is multiplicative in $n$, hence
$$
g(q,a)=\prod_{p}(1-\frac{(p^{k},q)}{p^{k}}\psi_{k}(p;q,a)).
$$
Let \bi f(q)=\prod_{p\mid q} (1-\frac{2}{p^{k}})^{-1},\ \
h(q,a)=\prod_{p\mid q}(1-\frac{(p^{k},q)}{p^{k}}\psi_{k}(p;q,a)),
\ei
then, by (3.8) of [2],
 \bi g(q,a)=\varrho f(q)h(q,a), \ei
where
\bi \varrho= \prod_{p} (1-\frac{2}{p^{k}}).\ei

By (3.9) and (3.10) of [2],  we have
 \bi h(q,a)=\prod_{\stackrel{p\mid
q}{(p^{k},q)\mid a(a+1)}}(1-\frac{(p^{k},q)}{p^{k}}),\ei
 and for co-prime natural numbers $q_{1}, q_{2}$
 \bi h(q_{1}q_{2},a)=h(q_{1},a)h(q_{2},a), \ei
  $h(q,a)$ is a  multiplicative function of $q$.\\
We need the following sums of Gaussian type
\bi
G(q,a)=\sum\limits_{b=1}^{q}g(q,b)e(\frac{ab}{q}),\
H(q,a)=\sum\limits_{b=1}^{q}h(q,b)e(\frac{ab}{q}),\
\ei
which by
(4.16) are related by
 \bi G(q,a)=\varrho f(q)H(q,a). \ei
 As in [2], we introduce the following function
 \bi H(q)= \sum\limits_{\stackrel{a=1}{(a,q)=1}}^{q}\mid
H(q,a)\mid^{2}.\ei
 We define
 \bi I(\beta)=\sum\limits_{n\leq x}e(\beta n),\ \
S^{\ast}(\alpha)=S^{\ast}(\alpha;q,a)=q^{-1}G(q,a)I(\alpha-\frac{a}{q}),
\ei
 \bi \Delta (\alpha)=\Delta (\alpha
;q,a)=S(\alpha)-S^{\ast}(\alpha).\ei
 Let
 \bi \mathfrak{M}(q,a;T)=\{ \alpha:\ \mid q\alpha-a\mid\leq T/x \},
\ei
and $ \mathfrak{M}(T)$ denote the union of the intervals
 $\mathfrak{M}(q,a;T)$ with $1\leq a\leq q\leq T$ and $(q,a)=1$, then
$$
 \mathfrak{M}(2T)\backslash  \mathfrak{M}(T)=\bigcup\limits_{\stackrel{T<q\leq 2T}{a\leq q,(a,q)=1}}
 \{\alpha :\ \mid q\alpha-a\mid\leq T/x\}
 \bigcup\limits_{\stackrel{q\leq T}{a\leq q,(a,q)=1}}\{\alpha :T/x< \mid q\alpha-a\mid\leq
 2T/x\}.
$$
and
\bi \mathfrak{M}=\mathfrak{M}(x/R)\subseteqq
\mathfrak{M}(1)\bigcup\limits_{1\leq 2^{j}\leq x/R}
\mathfrak{M}(2^{j+1})\backslash \mathfrak{M}(2^{j}). \ei
 Hence, as the proof of Lemma 4.2 in [1], \bi \int_{\mathfrak{M}}\mid
S^{\ast}(\alpha)\Delta (\alpha)\mid d\alpha \ll\log x
\max\limits_{1\leq T\leq
x/R}(U_{1}^{\ast}(T))^{1/2}(U_{2}^{\ast}(T))^{1/2}+\int_{\mathfrak{M}(1)}\mid
S^{\ast}(\alpha)\Delta (\alpha)\mid d\alpha, \ei
 where
  \bi
U_{1}^{\ast}(T)=\int_{\mathfrak{M}(2T)\backslash
\mathfrak{M}(T)}\mid S^{\ast}(\alpha)\mid^{2} d\alpha,\
U_{2}^{\ast}(T)=\int_{\mathfrak{M}(2T)}\mid \Delta (\alpha)\mid^{2}
d\alpha.\ei
 {\bf LEMMA 4.6.} Suppose that $1\leq T\leq \frac{1}{2}\sqrt{x}$.
Then
$$
\int_{\mathfrak{M}(T)}\mid S(\alpha)-S^{\ast}(\alpha)\mid^{2}d\alpha
\ll
T^{3-\frac{2}{k}}x^{\frac{2}{k}-1+\varepsilon}+x^{\frac{4}{k+1}-1+\varepsilon}T^{2}.
$$
Proof.\ This is Lemma 5.1 of [2].\\
{\bf LEMMA 4.7.}\ Suppose that $1\leq 2T\leq \frac{1}{2}\sqrt{x}$.
Then
$$
\int_{\mathfrak{M}(2T)\backslash  \mathfrak{M}(T)}\mid
S^{\ast}(\alpha)\mid^{2} d\alpha\ll xT^{\frac{1}{k}-1+\varepsilon}.
$$
Proof.\ This is Lemma 5.2 of [2].\\
{\bf LEMMA 4.8.} We have
$$
\int_{\mathfrak{M}}\mid S(\alpha)-S^{\ast}(\alpha)\mid^{2}d\alpha\ll
x^{2+\varepsilon}R^{\frac{2}{k}-3}+x^{\frac{4}{k+1}+1+\varepsilon}R^{-2}.
$$
Proof.\ By taking $T=x/R$ in Lemma 4.6, the lemma follows.\\
{\bf LEMMA 4.9.} We have
$$
\sum\limits_{q\leq
x/R}\sum\limits_{a=1}^{q}{}^{^{\prime}}\int_{\mathfrak{M}(q,a)}\mid
\Delta (\alpha ;q,a)\mid^{2}H_{q}(\alpha)d\alpha \ll
x^{\frac{5}{2}+\varepsilon}R^{\frac{2}{k}-3}+x^{\frac{3}{2}+\frac{4}{k+1}+\varepsilon}R^{-2}.
$$
Proof.\ By Lemma 4.5 and Lemma 4.6, the lemma follows.\\
{\bf LEMMA 4.10.}  We have
$$\int_{\mathfrak{M}(1)}\mid S^{\ast}(\alpha)\Delta (\alpha)\mid d\alpha
\ll x^{\frac{2}{k+1}+\varepsilon}.$$
Proof.\ By Cauchy-Schwarz inequality and Lemma 4.6 \\
$\int_{\mathfrak{M}(1)}\mid S^{\ast}(\alpha)\Delta (\alpha)\mid
d\alpha \ll\int_{1-1/x}^{1+1/x}\mid S^{\ast}(\alpha)\Delta
(\alpha)\mid d\alpha \\
\ll \int_{0}^{x^{-1}}\mid S^{\ast}(\alpha)\Delta (\alpha)\mid
d\alpha \ll
 (\int_{0}^{x^{-1}}\mid S^{\ast}(\alpha)\mid^{2}
d\alpha)^{1/2}(\int_{0}^{x^{-1}}\mid \Delta (\alpha)\mid^{2}
d\alpha)^{1/2}\\
 \ll
x^{\frac{1}{2}}(x^{\frac{2}{k}-1+\varepsilon}+x^{\frac{4}{k+1}-1+\varepsilon})^{1/2}
\ll x^{\frac{2}{k+1}+\varepsilon},$\\
the lemma follows.\\
By (4.28), Lemma 4.6 and Lemma 4.7, \bi U_{1}^{\ast}(T)\ll
xT^{\frac{1}{k}-1+\varepsilon},\ U_{2}^{\ast}(T)\ll
 T^{3-\frac{2}{k}}x^{\frac{2}{k}-1+\varepsilon}+x^{\frac{4}{k+1}-1+\varepsilon}T^{2}.\ei
 Hence, by (4.26), (4.27), (4.29) and Lemma 4.10,\\
 $ \int_{\mathfrak{M}}\mid S^{\ast}(\alpha)\Delta
(\alpha)\mid d\alpha \ll\log x \max\limits_{1\leq T\leq
x/R}(xT^{\frac{1}{k}-1+\varepsilon})^{1/2}
(T^{3-\frac{2}{k}}x^{\frac{2}{k}-1+\varepsilon}+x^{\frac{4}{k+1}-1+\varepsilon}T^{2})^{1/2}+
x^{\frac{2}{k+1}+\varepsilon}\\
\ll  x^{1+\frac{1}{2k}+\varepsilon}R^{\frac{1}{2k}-1}+
x^{\frac{2}{k+1}+\frac{1}{2k}+\frac{1}{2}+\varepsilon}R^{-\frac{1}{2k}-\frac{1}{2}}+x^{\frac{2}{k+1}+\varepsilon}.
$\\
Combining these with Lemma 4.5,\\
$ \sum\limits_{q\leq
x/R}\sum\limits_{a=1}^{q}{}^{^{\prime}}\int_{\mathfrak{M}(q,a)} \mid
S^{\ast}(\alpha;q,a)\Delta (\alpha ;q,a)H_{q}(\alpha)\mid
 d\alpha$\\
 \bi
 \ll x^{\frac{3}{2}+\frac{1}{2k}+\varepsilon}R^{\frac{1}{2k}-1}
 +x^{1+\frac{2}{k+1}+\frac{1}{2k}+\varepsilon}R^{-\frac{1}{2k}-\frac{1}{2}}+x^{\frac{1}{2}+\frac{2}{k+1}+\varepsilon}.
  \ei
 {\bf LEMMA 4.11.}\  We have
$$ \sum\limits_{q\leq
x/R}\sum\limits_{a=1}^{q}{}^{^{\prime}}\int_{\mathfrak{M}(q,a)} \mid
\Delta (\alpha ;q,a)\mid^{2}\mid K(\alpha)\mid d\alpha\ll
x^{1+\frac{1}{k}-\frac{\tau}{2}}+x^{\frac{1}{2}+\frac{4}{k+1}
-\frac{\tau}{2}}.$$\\
Proof.\  We use the method of (4.8) in [9]. Let
$\mathfrak{M}_{0}(1,1)=\{\alpha : \mid \alpha -1\mid\leq 1/x\}$
and  for $j\geq 1$\\
\[
\mathfrak{M}_{j}(q,a)=\left\{\begin{array}{cc} \{\alpha :
2^{j-1}x^{-1}<\mid
q\alpha -a\mid\leq 2^{j}x^{-1}\},& \ {\rm  when} \ \  q\leq 2^{j-1};\\
\{\alpha : \mid q\alpha -a\mid\leq  2^{j}x^{-1}\},& \ {\rm  when}\ \
2^{j-1}<q\leq 2^{j}.
\end{array}
\right.
\]
 $$ \mathfrak{M}_{j}=\bigcup_{\stackrel{1\leq a\leq q\leq
2^{j}}{(a,q)=1}}\mathfrak{M}_{j}(q,a). $$
Choose $J$ so that
$2^{J-1}<x/R\leq 2^{J}$. For $1\leq a\leq q \leq
x/R, \ (a,q)=1,$\\
then\\
$\mathfrak{M}(q,a)\subseteq\bigcup\limits_{\stackrel{0\leq j\leq
J}{q\leq 2^{j}}}\mathfrak{M}_{j}(q,a),\ \
\mathfrak{M}\subseteq \bigcup\limits_{0\leq j\leq J}\mathfrak{M}_{j}$.\\
By Lemme 4.4, when $\alpha\in \mathfrak{M}_{j}(q,a)$ with $1\leq
a\leq q\leq 2^{j}, (a,q)=1$,
 \bi K_{q}(\alpha)\ll 2^{-j}x\log x .\ei
By Lemma 4.6 and (4.31),\\
$ \sum\limits_{q\leq
x/R}\sum\limits_{a=1}^{q}{}^{^{\prime}}\int_{\mathfrak{M}(q,a)} \mid
\Delta (\alpha ;q,a)\mid^{2}\mid K_{q}(\alpha)\mid d\alpha \ll
 x^{\frac{2}{k}+\varepsilon}(\frac{x}{R})^{2-\frac{2}{k}}+x^{\frac{4}{k+1}+\varepsilon}(\frac{x}{R})$\\
 \bi
 \ll
 x^{2+\varepsilon}R^{\frac{2}{k}-2}+x^{1+\frac{4}{k+1}+\varepsilon}R^{-1}.
 \ei
 By Lemma 4.2, Lemma 4.9 and (4.32), the lemma follows.\\
By (4.23), (4.24) and Lemma 4.2, when $\alpha\in\mathfrak{M}(q,a)$,\\
$ \mid S(\alpha)\mid^{2}K(\alpha)=\mid
S^{\ast}(\alpha;q,a)\mid^{2}K(\alpha)\\
+2(\Re
(\overline{S^{\ast}(\alpha;q,a)}\Delta(\alpha;q,a)))(K_{q}(\alpha)+H_{q}(\alpha))+\mid\Delta(\alpha;q,a)\mid^{2}K(\alpha).
$\\
By (4.13), (4.30) and Lemma 4.11,
\bi
S_{A}(x,Q)=S_{A1}(x,Q)+S_{A2}(x,Q)+O(x^{1+\frac{1}{k}-\frac{\tau}{2}}+x^{\frac{1}{2}+\frac{4}{k+1}
-\frac{\tau}{2}}), \ei
 where
 \bi S_{A1}(x,Q)=\sum\limits_{q\leq
x/R}\sum\limits_{a=1}^{q}{}^{^{\prime}}\int_{\mathfrak{M}(q,a)} \mid
S^{\ast}(\alpha;q,a)\mid^{2}K(\alpha)d\alpha,
 \ei
$$ S_{A2}(x,Q)=2\sum\limits_{q\leq
x/R}\sum\limits_{a=1}^{q}{}^{^{\prime}}\int_{\mathfrak{M}(q,a)}(\Re
(\overline{S^{\ast}(\alpha;q,a)}\Delta(\alpha;q,a)))K_{q}(\alpha)d\alpha,
$$
by (4.21) and (4.23)
\bi S_{A2}(x,Q)=2\varrho\sum\limits_{q\leq
x/R}q^{-1}f(q)\sum\limits_{a=1}^{q}{}^{^{\prime}}\int_{-1/qR}^{1/qR}(\Re
(\overline{H(q,a)I(\beta)}\Delta(\frac{a}{q}+\beta;q,a)))K_{q}(\beta)d\beta.
\ei
\section{ The estimate of $S_{A2}(x,Q)$ }
\setcounter{equation}{0}
{\bf LEMMA 5.1.}\  We have
$$
\sum\limits_{b=1}^{q}h(q,b)\sum\limits_{r\mid
q}\mu(\frac{q}{r})g(r,b)
=q^{-1}\sum\limits_{b=1}^{q}h(q,b)\sum\limits_{a=1}^{q}{}^{\prime}G(q,a)e(\frac{-ab}{q}).
$$
Proof.\ By (4.20)
and $\sum\limits_{a=1}^{q}{}^{\prime}e(\frac{an}{q})=\sum\limits_{r\mid (q,n)}r\mu(\frac{q}{r})$, \\
$
q^{-1}\sum\limits_{b=1}^{q}h(q,b)\sum\limits_{a=1}^{q}{}^{\prime}G(q,a)e(\frac{-ab}{q})\\
=q^{-1}\sum\limits_{b=1}^{q}h(q,b)\sum\limits_{a=1}^{q}{}^{\prime}e(\frac{-ab}{q})
\sum\limits_{t=1}^{q}g(q,t)e(\frac{at}{q})\\
= q^{-1}\sum\limits_{b=1}^{q}h(q,b)\sum\limits_{r\mid
q}r\mu(\frac{q}{r}) \sum\limits_{\stackrel{t=1}{r\mid
t-b}}^{q}g(q,t).$\\
Hence \bi
q^{-1}\sum\limits_{b=1}^{q}h(q,b)\sum\limits_{a=1}^{q}{}^{\prime}G(q,a)e(\frac{-ab}{q})=q^{-1}\sum\limits_{r\mid
q}r\mu(\frac{q}{r})\sum\limits_{b=1}^{q}h(q,b)
\sum\limits_{\stackrel{t=1}{r\mid t-b}}^{q}g(q,t).\ei
 By (4.14),
\bi
 \sum\limits_{\stackrel{t=1}{t\equiv b(r)}}^{q}g(q,t)
 =\sum\limits_{n=1}^{\infty}\mu(n)\frac{(n^{k},q)}{n^{k}}
 \sum\limits_{\stackrel{t=1}{t\equiv b(r)}}^{q}\psi_{k}(n;q,t),\ei
where
$$ \psi_{k}(n;q,t)=\#\{u,v\geq 1: uv=n,\ (u,v)=1,\ (u^{k},q)\mid t,\
(v^{k},q)\mid t+1\ \}. $$
Write $t=b+rl, q=rm$, and $l$ runs through
complete residues
modulus $m$, we want to count the number of $l$.\\
We have\\
$ (u^{k},rm)\mid b+rl,\ (v^{k},rm)\mid b+1+rl . $\\
Let \\
$$
(u^{k},r)=r_{1},\ (v^{k},r)=r_{2}, (r_{1},r_{2})=1,$$
 then\\
$ r_{1}\mid b,\ r_{2}\mid b+1,  \
\\
(\frac{u^{k}}{r_{1}},m)\mid \frac{b}{r_{1}}+\frac{r}{r_{1}}l,\
(\frac{v^{k}}{r_{2}},m)\mid \frac{b+1}{r_{2}}+\frac{r}{r_{2}}l.$\\
These congruences have a unique solution of $l$  modulus
$(\frac{u^{k}}{r_{1}},m)(\frac{v^{k}}{r_{2}},m)$, therefore the
number of $l$ is
$\frac{m}{(\frac{u^{k}}{r_{1}},m)(\frac{v^{k}}{r_{2}},m)}=\frac{mr_{1}r_{2}}{(u^{k},r_{1}m)(v^{k},r_{2}m)}.$\\
 Since $(r_{1},r_{2})=(r_{1},v)=(r_{2},u)=1,\  $
the number of $l$ is\\
$
\frac{mr_{1}r_{2}}{(u^{k},r_{1}r_{2}m)(v^{k},r_{1}r_{2}m)}=\frac{mr_{1}r_{2}}{(n^{k},r_{1}r_{2}m)}.$\\
Hence\\
$ \sum\limits_{\stackrel{t=1}{t\equiv b(r)}}^{q}g(q,t)
=\sum\limits_{n=1}^{\infty}\mu(n)\frac{(n^{k},q)}{n^{k}}\sum\limits_{u,v}{}^{\ast}\frac{mr_{1}r_{2}}{(n^{k},r_{1}r_{2}m)},
$\\
where  $\sum\limits_{u,v}{}^{\ast}$ means $(u^{k},r)=r_{1}\mid b, \
(v^{k},r)=r_{2}\mid b+1,\ uv=n, (u,v)=1,\ u,v\geq 1, q=rm.$\\
Let $r=r_{1}r_{2}\Delta,$ \ then \\
$(q,n^{k})=(r_{1}r_{2}\Delta m,u^{k}v^{k})=(r_{1}r_{2}\Delta
m,u^{k})(r_{1}r_{2}\Delta m,v^{k}) =(r_{1}\Delta
m,u^{k})(r_{2}\Delta m,v^{k})\\
=r_{1}r_{2}(\Delta
m,\frac{u^{k}}{r_{1}})(\Delta m,\frac{v^{k}}{r_{2}}).$\\
Since
$(\frac{u^{k}}{r_{1}},r_{2}\Delta)=(\frac{v^{k}}{r_{2}},r_{1}\Delta)=1$,
 we have
$(\Delta,\frac{u^{k}}{r_{1}}\frac{v^{k}}{r_{2}})=1,$\\
hence\\
$ (q,n^{k})=r_{1}r_{2}(
m,\frac{u^{k}}{r_{1}})(m,\frac{v^{k}}{r_{2}})=(r_{1}r_{2}m,u^{k})(r_{1}r_{2}m,v^{k})=(r_{1}r_{2}m,(uv)^{k}).$\\
So, by (5.1) and (5.2), \\
$
q^{-1}\sum\limits_{b=1}^{q}h(q,b)\sum\limits_{a=1}^{q}{}^{\prime}G(q,a)e(\frac{-ab}{q})
=\sum\limits_{r\mid q}\mu(\frac{q}{r})\sum\limits_{b=1}^{q}h(q,b)
\sum\limits_{n=1}^{\infty}\mu(n)\frac{1}{n^{k}}\sum\limits_{u,v}{}^{\ast}\frac{rmr_{1}r_{2}}{q}\\
=\sum\limits_{r\mid q}\mu(\frac{q}{r})\sum\limits_{b=1}^{q}h(q,b)
\sum\limits_{n=1}^{\infty}\mu(n)\frac{1}{n^{k}}\sum\limits_{u,v}{}^{\ast}r_{1}r_{2}
=\sum\limits_{r\mid q}\mu(\frac{q}{r})\sum\limits_{b=1}^{q}h(q,b)
\sum\limits_{n=1}^{\infty}\mu(n)\frac{(n^{k},r)}{n^{k}}\psi_{k}(n;r,b),
$\\
by (4.14), the lemma follows.\\
Write
$$L=\sum\limits_{a=1}^{q}{}^{^{\prime}}\overline{H(q,a)}\Delta(\frac{a}{q}+\beta;q,a),$$
then, by (4.20) and (4.24),\\
$L=\sum\limits_{b=1}^{q}h(q,b)\sum\limits_{a=1}^{q}{}^{^{\prime}}e(\frac{-ab}{q})\Delta(\frac{a}{q}+\beta;q,a)\\
=\sum\limits_{n\leq
x}(\mu_{k}(n)\mu_{k}(n+1)\sum\limits_{b=1}^{q}h(q,b)\sum\limits_{a=1}^{q}{}^{^{\prime}}e(\frac{a}{q}(n-b))-
q^{-1}\sum\limits_{b=1}^{q}h(q,b)\sum\limits_{a=1}^{q}{}^{^{\prime}}G(q,a)e(\frac{-ab}{q}))e(\beta
n).$\\
Write\\
$U(y)=\sum\limits_{n\leq
y}(\mu_{k}(n)\mu_{k}(n+1)\sum\limits_{b=1}^{q}h(q,b)\sum\limits_{a=1}^{q}{}^{^{\prime}}e(\frac{a}{q}(n-b))-
q^{-1}\sum\limits_{b=1}^{q}h(q,b)\sum\limits_{a=1}^{q}{}^{^{\prime}}G(q,a)e(\frac{-ab}{q}))\\
=U_{1}(y)-U_{2}(y). $\\
and using $e(n\beta)=e(x\beta)-2\pi i\beta\int_{n}^{x}e(y\beta)dy$,\
we have \bi L=U(x)e(x\beta)-2\pi
i\beta\int_{1}^{x}U(y)e(y\beta)dy.
\ei
 We estimate $U_{1}(y)$ first.\\
Write
 \bi
 E(x;q,a)= A_{k}(x;q,a)-q^{-1}g(q,a)x. \ei
We have\\
 $ U_{1}(y)=U_{1}([y])=\sum\limits_{b=1}^{q}h(q,b)\sum\limits_{n\leq
[y]}\mu_{k}(n)\mu_{k}(n+1)\sum\limits_{r\mid
(q,n-b)}r\mu(\frac{q}{r})\\
=\sum\limits_{r\mid q}r\mu(\frac{q}{r})
\sum\limits_{b=1}^{q}h(q,b)\sum\limits_{\stackrel{n\leq [y]}{r\mid
n-b}}\mu_{k}(n)\mu_{k}(n+1)\\
=\sum\limits_{r\mid q}\mu(\frac{q}{r})
\sum\limits_{b=1}^{q}h(q,b)g(r,b)[y]+\sum\limits_{r\mid
q}r\mu(\frac{q}{r}) \sum\limits_{b=1}^{q}h(q,b)E([y];r,b),$\\
and\\
$ U_{2}(y)= \sum\limits_{n\leq [y]}
q^{-1}\sum\limits_{b=1}^{q}h(q,b)\sum\limits_{a=1}^{q}{}^{^{\prime}}G(q,a)e(\frac{-ab}{q})\\
=q^{-1}\sum\limits_{b=1}^{q}h(q,b)\sum\limits_{a=1}^{q}{}^{^{\prime}}G(q,a)e(\frac{-ab}{q})[y].
$\\
We need an obvious result on $Y_{k}(x,Q)$.\\
{\bf LEMMA 5.2.}\ For all  $Q,x$, we have
$$Y_{k}(x,Q)\ll x^{2}\log Q +Q^{2}.$$
Proof.\ By (1.4) and $g(q,a)\ll 1$, \\
$Y_{k}(x,Q)\ll \sum\limits_{q\leq
Q}\sum\limits_{a=1}^{q}(x/q+1)^{2}\ll \sum\limits_{q\leq
Q}\sum\limits_{a=1}^{q}(x^{2}q^{-2}+1)\ll x^{2}\log Q +Q^{2},$\\
the lemma follows.\\
By Lemma 5.1, we have\\
$U(y)=\sum\limits_{r\mid
q}r\mu(\frac{q}{r}) \sum\limits_{b=1}^{q}h(q,b)E([y];r,b). $\\
Hence\\
 $ L\ll \sum\limits_{r\mid q}r\mu^{2}(\frac{q}{r})
\sum\limits_{b=1}^{q}h(q,b)\mid E([x];r,b)\mid + \mid \beta\mid
\int_{1}^{x}\sum\limits_{r\mid q}r\mu^{2}(\frac{q}{r})
\sum\limits_{b=1}^{q}h(q,b)\mid E([y];r,b)\mid dy. $\\
By (4.35), Lemma 4.4 and $I(\beta)\ll x/(1+x\mid \beta\mid)$\\
$S_{A2}(x,Q)\ll \sum\limits_{q\leq
x/R}q^{-1}\int_{0}^{1/qR}\frac{x}{(1+x\beta)}\frac{x\log x}{q(1+x
\beta)}(\sum\limits_{r\mid q}r\mu^{2}(\frac{q}{r})
\sum\limits_{b=1}^{q}h(q,b)\mid E([x];r,b)\mid \\
+  \beta \int_{1}^{x}\sum\limits_{r\mid q}r\mu^{2}(\frac{q}{r})
\sum\limits_{b=1}^{q}h(q,b)\mid E([y];r,b)\mid dy)d\beta  \ll
\int_{0}^{R^{-1}}\frac{x^{2}\log x}{(1+x
\beta)^{2}}M(\beta)d\beta,$\\
where \\
$ M(\beta)\ll$ \bi \sum\limits_{\stackrel{q\leq x/R}{q\leq 1/\beta
R}}q^{-2}\\
(\sum\limits_{r\mid q}r\mu^{2}(\frac{q}{r})
\sum\limits_{b=1}^{q}h(q,b)\mid E([x];r,b)\mid +  \beta
\int_{1}^{x}\sum\limits_{r\mid q}r\mu^{2}(\frac{q}{r})
\sum\limits_{b=1}^{q}h(q,b)\mid E([y];r,b)\mid dy).\ei
 When $r\mid q$, we have $\sum\limits_{b=1}^{q}\mid E(y;r,b)\mid=qr^{-1}\sum\limits_{b=1}^{r}\mid E(y;r,b)\mid$, then\\
$ \sum\limits_{\stackrel{q\leq x/R}{q\leq 1/\beta
R}}q^{-2}\sum\limits_{r\mid q}r\mu^{2}(\frac{q}{r})
\sum\limits_{b=1}^{q}h(q,b)\mid E([x];r,b)\mid \ll
\sum\limits_{r\leq x/R}\sum\limits_{b=1}^{r}\mid E([x];r,b)\mid
\sum\limits_{\stackrel{q\leq x/R}{r\mid q}}q^{-1}\\
\ll \log x\sum\limits_{r\leq x/R}r^{-1}\sum\limits_{b=1}^{r}\mid
E([x];r,b)\mid.$\\
By Cauchy inequality and (1.5),\\
$ \sum\limits_{\stackrel{q\leq x/R}{q\leq 1/\beta
R}}q^{-2}\sum\limits_{r\mid q}r\mu^{2}(\frac{q}{r})
\sum\limits_{b=1}^{q}h(q,b)\mid E([x];r,b)\mid \\
\ll \log x (\sum\limits_{r\leq
x/R}\sum\limits_{b=1}^{r}r^{-2})^{1/2}(\sum\limits_{r\leq
x/R}\sum\limits_{b=1}^{r}\mid E([x];r,b)\mid^{2})^{1/2}\\
\ll
\log^{3/2}x(x^{\frac{2}{k}+\varepsilon}(\frac{x}{R})^{2-\frac{2}{k}}+x^{\frac{4}{k+1}+\varepsilon}
 )^{1/2}
 \ll
 x^{1+\varepsilon}R^{\frac{1}{k}-1}+x^{\frac{2}{k+1}+\varepsilon},$\\
this contributes to $S_{A2}(x,Q)$\\
$\ll
(x^{1+\varepsilon}R^{\frac{1}{k}-1}+x^{\frac{2}{k+1}+\varepsilon})\int_{0}^{R^{-1}}\frac{x^{2}\log
x}{(1+x\beta)^{2}}d\beta \\
\ll x\log
x(x^{1+\varepsilon}R^{\frac{1}{k}-1}+x^{\frac{2}{k+1}+\varepsilon})
\ll
x^{2+\varepsilon}R^{\frac{1}{k}-1}+x^{1+\frac{2}{k+1}+\varepsilon}.$\\
For the second term of $M(\beta)$, we have\\
$\sum\limits_{q\leq x/R}q^{-2}\beta \int_{1}^{x}\sum\limits_{r\mid
q}r\mu^{2}(\frac{q}{r}) \sum\limits_{b=1}^{q}h(q,b)\mid
E([y];r,b)\mid dy\\
 \ll \beta\int_{1}^{x} \log
x\sum\limits_{r\leq x/R}r^{-1}\sum\limits_{b=1}^{r}\mid
E([y];r,b)\mid dy\\
\ll  \beta\int_{1}^{x/R} \log x\sum\limits_{r\leq
x/R}r^{-1}\sum\limits_{b=1}^{r}\mid E([y];r,b)\mid dy +
\beta\int_{x/R}^{x} \log x\sum\limits_{r\leq
x/R}r^{-1}\sum\limits_{b=1}^{r}\mid E([y];r,b)\mid dy .$\\
By Cauchy inequality and Lemma 5.2,\\
$\beta\int_{1}^{x/R} \log^{3/2} x(\sum\limits_{r\leq
x/R}\sum\limits_{b=1}^{r}\mid E([y];r,b)\mid^{2})^{1/2} dy\\
\ll \beta\int_{1}^{x/R}\log^{3/2}x(y^{2}\log x+(\frac{x}{R})^{2})^{1/2}dy \\
\ll x^{2}R^{-2}\beta\log^{2} x, $\\
this contributes to $S_{A2}(x,Q)$\\
$\ll x^{2}R^{-2} \log^{2} x\int_{0}^{R^{-1}}\frac{x^{2}\log
x}{(1+x\beta)^{2}} \beta d\beta
\ll x^{2}R^{-2} \log^{4} x.$\\
By Cauchy inequality and (1.5),\\
 $
\beta\int_{x/R}^{x} \log x\sum\limits_{r\leq
x/R}r^{-1}\sum\limits_{b=1}^{r}\mid E([y];r,b)\mid dy \\
\ll \beta\int_{x/R}^{x}\log^{3/2}x(x^{\frac{2}{k}+\varepsilon}
(\frac{x}{R})^{2-\frac{2}{k}}+x^{\frac{4}{k+1}+\varepsilon})^{1/2}dy
\ll \beta(
x^{2+\varepsilon}R^{\frac{1}{k}-1}+x^{1+\frac{2}{k+1}+\varepsilon}),$\\
this contributes to $S_{A2}(x,Q)$\\
$$ \ll
x^{2+\varepsilon}R^{\frac{1}{k}-1}+x^{1+\frac{2}{k+1}+\varepsilon}.$$
 Consequently,
\bi S_{A2}(x,Q)\ll
x^{2+\varepsilon}R^{\frac{1}{k}-1}+x^{1+\frac{2}{k+1}+\varepsilon}+x^{2}R^{-2}
\log^{4} x.\ei
\section{ The estimate of $S_{A1}(x,Q)$ }
\setcounter{equation}{0}
By (4.2), (4.23) and (4.34),\\
$ S_{A1}(x,Q)=\sum\limits_{q\leq
x/R}\sum\limits_{a=1}^{q}{}^{^{\prime}}\int_{-1/2}^{1/2}q^{-2} \mid
G(q,a)I(\beta)\mid^{2}K(\frac{a}{q}+\beta)d\beta -\\
\sum\limits_{q\leq
x/R}\sum\limits_{a=1}^{q}{}^{^{\prime}}\int_{1/qR\leq\mid
\beta\mid\leq 1/2}q^{-2} \mid
G(q,a)I(\beta)\mid^{2}K(\frac{a}{q}+\beta)d\beta\\ $ \bi
=S_{B}(x,Q)-S_{C}(x,Q),
 \ei
say.\\
We consider $S_{C}(x,Q)$ first. By Lemma 4.4,  \\
 $ S_{C}(x,Q)\ll   \sum\limits_{q\leq
x/R}q^{-2}\sum\limits_{a=1}^{q}{}^{^{\prime}}\mid
G(q,a)\mid^{2}\int_{1/qR\leq\mid \beta\mid\leq
1/2}\min(x^{2},\beta^{-2})\sum\limits_{u\leq
\sqrt{x}}\frac{x}{u+x\parallel u\frac{a}{q}+u\beta\parallel}
d\beta\\
\ll  \sum\limits_{q\leq
x/R}q^{-2}\sum\limits_{a=1}^{q}{}^{^{\prime}}\mid
G(q,a)\mid^{2}\sum\limits_{u\leq \sqrt{x}}\int_{1/qR\leq\mid
\beta\mid\leq 1/2}\frac{x\beta^{-2}}{u+x\parallel
u\frac{a}{q}+u\beta\parallel}
d\beta\\
\ll\sum\limits_{q\leq
x/R}q^{-2}\sum\limits_{a=1}^{q}{}^{^{\prime}}\mid
G(q,a)\mid^{2}\sum\limits_{u\leq \sqrt{x}}\int_{u/qR\leq\mid
\gamma\mid\leq u/2}u\frac{x\gamma^{-2}}{u+x\parallel
u\frac{a}{q}+\gamma\parallel} d\gamma. $\\
Now,\\
$ux\sum\limits_{0<\mid j\mid\leq u/2}\int_{j-1/2\leq \gamma\leq
j+1/2}\frac{\gamma^{-2}}{u+x\parallel u\frac{a}{q}+\gamma\parallel}
d\gamma
 \ll
 ux\sum\limits_{0<\mid j\mid\leq u/2}\int_{-1/2\leq
\gamma\leq 1/2}\frac{j^{-2}}{u+x\parallel
u\frac{a}{q}+\beta+j\parallel} d\beta\\
\ll ux\sum\limits_{0<\mid j\mid\leq u/2}j^{-2}\int_{-1/2\leq
\gamma\leq 1/2}\frac{1}{u+x\mid \beta\mid} d\beta
\ll u\log x.$\\
Hence \\
$ S_{C}(x,Q)\ll \sum\limits_{q\leq
x/R}q^{-2}\sum\limits_{a=1}^{q}{}^{^{\prime}}\mid
G(q,a)\mid^{2}\sum\limits_{u\leq \sqrt{x}}u\log x\\
+\sum\limits_{q\leq
x/R}q^{-2}\sum\limits_{a=1}^{q}{}^{^{\prime}}\mid
G(q,a)\mid^{2}\sum\limits_{u\leq \sqrt{x}}\int_{u/qR\leq\mid
\gamma\mid\leq 1/2}u\frac{x\gamma^{-2}}{u+x\parallel
u\frac{a}{q}+\gamma\parallel} d\gamma\\
=S_{C1}(x,Q)+ S_{C2}(x,Q), $\ \   say.\\
By Lemma 4.2 of [2], for any $Q\geq 1$,\\
\bi
\sum\limits_{Q<q\leq 2Q}q^{-2}f^{2}(q)H(q)\ll
Q^{\frac{1}{k}-1+\varepsilon}. \ei Hence \bi S_{C1}(x,Q)\ll x\log
x\sum\limits_{q\leq x/R}q^{-2}f^{2}(q)H(q)\ll x\log x.
 \ei
{\bf LEMMA 6.1.}\ We have
$$S_{C2}(x,Q)\ll x^{\frac{3}{2}+\frac{1}{2k}+2\tau}.$$
Proof.\ We have\\
$S_{C2}(x,Q)\ll \sum\limits_{q\leq
x/R}q^{-2}\sum\limits_{a=1}^{q}{}^{^{\prime}}\mid
G(q,a)\mid^{2}\sum\limits_{u\leq \sqrt{x}}x\int_{u/qR\leq\mid
\gamma\mid\leq 1/2}\gamma^{-2}d\gamma\\\
\ll \sum\limits_{q\leq
x/R}q^{-2}\sum\limits_{a=1}^{q}{}^{^{\prime}}\mid
G(q,a)\mid^{2}\sum\limits_{u\leq \sqrt{x}}x\frac{qR}{u}\\
\ll xR\log x\sum\limits_{q\leq
x/R}q^{-1}\sum\limits_{a=1}^{q}{}^{^{\prime}}\mid G(q,a)\mid^{2}.$\\
 By partial summation and (6.2),\\
 $S_{C2}(x,Q)\ll xR\log x\sum\limits_{q\leq
x/R}q^{-1}H(q)f(q)^{2}\ll xR(\frac{x}{R})^{\frac{1}{k}+\varepsilon}
\ll x^{\frac{3}{2}+\frac{1}{2k}+2\tau},$\\
the lemma follows.\\
By Lemma 6.1 and (6.3),
$$ S_{C}(x,Q)\ll x^{\frac{3}{2}+\frac{1}{2k}+2\tau}.$$
Hence
\bi
 S_{A1}(x,Q)=S_{B}(x,Q)+O( x^{\frac{3}{2}+\frac{1}{2k}+2\tau}).
\ei
\section{ The estimate of $S_{B}(x,Q)$ }
\setcounter{equation}{0}
We have\\
$ \int_{-1/2}^{1/2}\mid I(\beta)\mid^{2}K(\frac{a}{q}+\beta)d\beta
=\sum\limits_{n_{1},n_{2}\leq x}\sum\limits_{\stackrel{uv\leq
x,u\leq Q}{n_{1}-n_{2}+uv=0}}e(uv\frac{a}{q}).$\\
By (6.1), \\
$S_{B}(x,Q)=\sum\limits_{q\leq
x/R}q^{-2}\sum\limits_{a=1}^{q}{}^{^{\prime}}\mid
G(q,a)\mid^{2}\sum\limits_{u\leq Q}\sum\limits_{v\leq
x/u}e(uv\frac{a}{q})\sum\limits_{\stackrel{m,n\leq x}{n=m+uv}}1. $\\
Hence
\bi S_{B}(x,Q)=\sum\limits_{q\leq
x/R}q^{-2}\sum\limits_{a=1}^{q}{}^{^{\prime}}\mid
G(q,a)\mid^{2}\sum\limits_{u\leq Q}\sum\limits_{v\leq
x/u}e(uv\frac{a}{q})[x-uv].
 \ei
By (6.2),\\
$
 \sum\limits_{u\leq Q}\sum\limits_{v\leq
x/u}(x-uv)\sum\limits_{q>
x/R}q^{-2}\sum\limits_{a=1}^{q}{}^{^{\prime}}\mid G(q,a)\mid^{2}\\
\ll \sum\limits_{u\leq Q}\sum\limits_{v\leq x/u}(x-uv)(
x/R)^{\frac{1}{k}-1+\varepsilon}\ll
(x/R)^{\frac{1}{k}-1+\varepsilon}\sum\limits_{u\leq
Q}u\sum\limits_{v\leq x/u}(x/u-v)\\
\ll (x/R)^{\frac{1}{k}-1+\varepsilon}\sum\limits_{u\leq Q}u(
x/u)^{2}
\ll  x^{\frac{3}{2}+\frac{1}{2k}+2\tau},$\\
and\\
$ \sum\limits_{q\leq
x/R}q^{-2}\sum\limits_{a=1}^{q}{}^{^{\prime}}\mid
G(q,a)\mid^{2}\sum\limits_{u\leq Q}\sum\limits_{v\leq
x/u}1\ll x\log x $.\\
Hence
 \bi
 S_{B}(x,Q)=\sum\limits_{u\leq Q}\sum\limits_{v\leq
x/u}(x-uv)\mathfrak{S}(uv)+O(x^{\frac{3}{2}+\frac{1}{2k}+2\tau})=S_{E}(x,Q)+O(x^{\frac{3}{2}+\frac{1}{2k}+2\tau}),
\ei
where
 \bi
 \mathfrak{S}(n)=\sum\limits_{q=1}^{\infty}\sum\limits_{a=1}^{q}{}^{\prime}q^{-2}\mid
 G(q,a)\mid^{2}e(\frac{an}{q}),
  \ei

\bi S_{E}(x,Q)=\sum\limits_{u\leq Q}\sum\limits_{v\leq
x/u}(x-uv)\mathfrak{S}(uv). \ei

\section{ The formula for $\mathfrak{S}(n)$ }
\setcounter{equation}{0}
 Write
\bi
 J(q)=J(q;n)=\sum\limits_{a=1}^{q}{}^{\prime}\mid
 H(q,a)\mid^{2}e(\frac{an}{q}),
 \ei
 then, by (4.21) and (7.3),
$$\mathfrak{S}(n)=\varrho^{2}\sum\limits_{q=1}^{\infty}q^{-2}f(q)^{2}J(q).$$
By (6.2),  $\mathfrak{S}(n)$ converges absolutely.\\
{\bf LEMMA 8.1.}\ For any natural numbers $q_{1}, q_{2},\ (q_{1},q_{2})=1$, we have\\
$$J(q_{1}q_{2})=J(q_{1})J(q_{2}),$$
the function $J(q)$ is multiplicative.\\
Proof.\ By (4.19) and $h(q,a+mq)=h(q,a)$, write $a=a_{1}q_{2}+a_{2}q_{1},\ b=b_{1}q_{2}+b_{2}q_{1}$,\\
$
J(q_{1}q_{2})\\
=\sum\limits_{a_{1}=1}^{q_{1}}{}^{\prime}\sum\limits_{a_{2}=1}^{q_{2}}{}^{\prime}
\mid\sum\limits_{b_{1}=1}^{q_{1}}{}^{\prime}h(q_{1},b_{1}q_{2})e(\frac{a_{1}q_{2}b_{1}}{q_{1}})\mid^{2}
e(\frac{a_{1}n}{q_{1}})\mid\sum\limits_{b_{2}=1}^{q_{2}}{}^{\prime}h(q_{2},b_{2}q_{1})e(\frac{a_{2}q_{1}b_{2}}{q_{2}})
\mid^{2} e(\frac{a_{2}n}{q_{2}})\\
=\sum\limits_{a_{1}=1}^{q_{1}}{}^{\prime}
\mid\sum\limits_{b_{1}=1}^{q_{1}}{}^{\prime}h(q_{1},b_{1}q_{2})e(\frac{a_{1}q_{2}b_{1}}{q_{1}})\mid^{2}
e(\frac{a_{1}n}{q_{1}})
 \sum\limits_{a_{2}=1}^{q_{2}}{}^{\prime}
\mid\sum\limits_{b_{2}=1}^{q_{2}}{}^{\prime}h(q_{2},b_{2}q_{1})e(\frac{a_{2}q_{1}b_{2}}{q_{2}})
\mid^{2} e(\frac{a_{2}n}{q_{2}})\\
=J(q_{1})J(q_{2}),$\\
the lemma follows.\\
We need to give a formula for $J(p^{t})$, here $p$ is a prime number and $t\geq 1$ .\\
We have\\
\bi J(p^{t})= \sum\limits_{a=1}^{p^{t}}\mid
 H(p^{t},a)\mid^{2}e(\frac{an}{p^{t}})
 -\sum\limits_{\stackrel{a=1}{p\mid a}}^{p^{t}}\mid
 H(p^{t},a)\mid^{2}e(\frac{an}{p^{t}})=J_{1}(p^{t})-J_{2}(p^{t}),
 \ei
where\\
$$J_{1}(p^{t})=\sum\limits_{a=1}^{p^{t}}\mid
 H(p^{t},a)\mid^{2}e(\frac{an}{p^{t}}),$$
$$J_{2}(p^{t})=\sum\limits_{\stackrel{a=1}{p\mid a}}^{p^{t}}\mid
 H(p^{t},a)\mid^{2}e(\frac{an}{p^{t}}).$$
{\bf LEMMA 8.2.}\ When $t>k$
$$ J(p^{t})=0.  $$
Proof.\ When  $t>k$, by (4.20),\\
$J_{1}(p^{t})=\sum\limits_{a=1}^{p^{t}}
 \mid\sum\limits_{b=1}^{p^{t}}h(p^{t},b)e(\frac{ab}{p^{t}})\mid^{2}e(\frac{an}{p^{t}})\\
=\sum\limits_{b_{1},b_{2}=1}^{p^{t}}h(p^{t},b_{1})h(p^{t},b_{2})\sum\limits_{a=1}^{p^{t}}
e(\frac{a(b_{1}-b_{2})}{p^{t}})e(\frac{an}{p^{t}})\\
=p^{t}\sum\limits_{\stackrel{b_{1},b_{2}=1}{b_{1}-b_{2}+n\equiv
0(p^{t})}}^{p^{t}}h(p^{t},b_{1})h(p^{t},b_{2})\\
=p^{t}\sum\limits_{b=1}^{p^{t}}h(p^{t},b)h(p^{t},b+n). $\\
Also\\
 $J_{2}(p^{t})=\sum\limits_{b_{1},b_{2}=1}^{p^{t}}h(p^{t},b_{1})h(p^{t},b_{2})\sum\limits_{\stackrel{a=1}{p\mid
a}}^{p^{t}} e(\frac{a(b_{1}-b_{2}+n)}{p^{t}})\\
=\sum\limits_{b_{1},b_{2}=1}^{p^{t}}h(p^{t},b_{1})h(p^{t},b_{2})
\sum\limits_{a=1}^{p^{t-1}} e(\frac{a(b_{1}-b_{2}+n)}{p^{t-1}})\\
=p^{t-1}\sum\limits_{\stackrel{b_{1},b_{2}=1}{b_{1}-b_{2}+n\equiv
0(p^{t-1})}}^{p^{t}}h(p^{t},b_{1})h(p^{t},b_{2}).$\\
By (3.9) of [2], for $t>k$, one has $h(p^{t},a)=h(p^{t},b)$ whenever
$a\equiv b(mod \ p^{k})$.Hence\\
$J_{2}(p^{t})
=p^{t-1}\sum\limits_{b=1}^{p^{t}}h(p^{t},b)h(p^{t},b+n)
\sum\limits_{\stackrel{b_{2}=1}{b_{2}\equiv b+n (p^{t-1})}}^{p^{t}}
=p^{t}\sum\limits_{b=1}^{p^{t}}h(p^{t},b)h(p^{t},b+n)=J_{1}(p^{t}),
$\\
therefore
$$J(p^{t})=J_{1}(p^{t})-J_{2}(p^{t})=0,$$
the lemma follows.\\
From now on, we suppose $1\leq t\leq k.$\\
By (4.7) of [2],
\bi
h(p^{t},a)=\left\{\begin{array}{cc}
1-p^{t-k},&\ {\rm \ if}\
p^{t}\mid a(a+1),\\
1,&\ {\rm \ otherwise}.
\end{array}
\right.
\ei
 For fixed $p^{t}, n$, write \\
$$\Phi_{1}=\Phi_{1}(p^{t})=\#\{b:\ p^{t}\mid b(b+1),\ b=1,2,...,p^{t} \}, $$
\bi \Phi_{2}=\Phi_{2}(p^{t})=\#\{b:\ p^{t}\mid b(b+1),\ p^{t}\mid
(b+n)(b+n+1),\ b=1,2,...,p^{t} \},\ei
$$ \Phi_{3}=\Phi_{3}(p^{t})=\#\{b:\ p^{t}\dagger b(b+1),\ p^{t}\mid (b+n)(b+n+1),\
b=1,2,...,p^{t} \},$$
$$ \Phi_{4}=\Phi_{4}(p^{t})=\#\{b:\ p^{t}\dagger b(b+1),\
b=1,2,...,p^{t} \},$$
\bi \Phi_{5}=\Phi_{5}(p^{t})=\#\{b_{1},
b_{2}:\ p^{t}\mid b_{1}(b_{1}+1),\ p^{t}\mid b_{2}(b_{2}+1),
b_{1}-b_{2}+n\equiv 0(p^{t-1}),\ b_{1},b_{2}=1,2,...,p^{t} \}.\ei
 {\bf LEMMA 8.3.}\ Suppose that $1\leq t\leq k$.  Then
$$J_{1}(p^{t})=p^{2t}-4p^{2t-k}+p^{3t-2k}\Phi_{2}.$$
Proof.\ As in the proof of Lemma 8.1 and Lemma 8.2, by (8.3),\\
$J_{1}(p^{t})=p^{t}\sum\limits_{b=1}^{p^{t}}h(p^{t},b)h(p^{t},b+n)
=p^{t}(J_{A}+J_{B}),$\\
where\\
$J_{A}=\sum\limits_{\stackrel{b=1}{p^{t}\mid
b(b+1)}}^{p^{t}}(1-p^{t-k})h(p^{t},b+n),\
J_{B}=\sum\limits_{\stackrel{b=1}{p^{t}\dagger
b(b+1)}}^{p^{t}}h(p^{t},b+n).$\\
Firstly, we deal with $J_{A}$.\\
We note that $\Phi_{1}=2,\ \Phi_{3}=2-\Phi_{2}, \Phi_{4}=p^{t}-2.$\\
By (8.3),\\
 $J_{A}=(1-p^{t-k})^{2}\Phi_{2}+(1-p^{t-k})\sum\limits_{\stackrel{b=1}{p^{t}\mid b(b+1),p^{t}\dagger (b+n)(b+n+1)}}^{p^{t}}1
 \\
 =(1-p^{t-k})^{2}\Phi_{2}+(1-p^{t-k})(2-\Phi_{2})
=2(1-p^{t-k})+(1-p^{t-k})(-p^{t-k})\Phi_{2}, $\\
hence
$$J_{A}=2(1-p^{t-k})-p^{t-k}(1-p^{t-k})\Phi_{2}.$$
Now, we deal with $J_{B}$. By (8.3),\\
$J_{B}=(1-p^{t-k})\Phi_{3}+\Phi_{4}-\Phi_{3}=p^{t}-2-p^{t-k}\Phi_{3}
=p^{t}-2-p^{t-k}(2-\Phi_{2}),$\\
hence\\
$$J_{B}=p^{t}-2-2p^{t-k}+p^{t-k}\Phi_{2}.$$
Therefore\\
$J_{1}(p^{t})=p^{t}(2(1-p^{t-k})-p^{t-k}(1-p^{t-k})\Phi_{2}+p^{t}-2-2p^{t-k}+p^{t-k}\Phi_{2})
=p^{2t}-4p^{2t-k}+p^{3t-2k}\Phi_{2},$\\
the lemma follows.\\
{\bf LEMMA 8.4.}\ We have
$$
J_{2}(p)=p^{2}(1-2p^{-k})^{2}.
$$
Proof.\ By (4.20) and (8.2), \\
$J_{2}(p)=(\sum\limits_{b=1}^{p}h(p,b))^{2},$\\
and\\
$ \sum\limits_{b=1}^{p}h(p,b)=\sum\limits_{\stackrel{b=1}{p\mid
b(b+1)}}^{p}(1-p^{1-k})+\sum\limits_{\stackrel{b=1}{p\dagger
b(b+1)}}^{p}1 =2(1-p^{1-k})+p-2=p-2p^{1-k}.$\\
Hence
$$J_{2}(p)=p^{2}(1-2p^{-k})^{2},$$
the lemma follows.\\
{\bf LEMMA 8.5.}\ Suppose that $1< t\leq k$. Then
$$J_{2}(p^{t})=p^{t-1}(p^{t+1}-4p^{t-k+1}+p^{2t-2k} \Phi_{5}).$$
Proof.\ By (4.20) and (8.2),\\
$J_{2}(p^{t})=p^{t-1}M,$\\
where
$$M=\sum\limits_{\stackrel{b_{1},b_{2}=1}{b_{1}-b_{2}+n\equiv
0(p^{t-1})}}^{p^{t}}h(p^{t},b_{1})h(p^{t},b_{2}).$$
By (8.3), we obtain\\
$M=M_{A}+M_{B},$\\
where\\
 $M_{A}=\sum\limits_{\stackrel{b_{1}=1}{p^{t}\mid b_{1}(b_{1}+1)}}^{p^{t}}
 \sum\limits_{\stackrel{b_{2}=1}{b_{1}-b_{2}+n\equiv
0(p^{t-1})}}^{p^{t}}(1-p^{t-k})h(p^{t},b_{2}),$\\
$M_{B}=\sum\limits_{\stackrel{b_{1}=1}{p^{t}\dagger
b_{1}(b_{1}+1)}}^{p^{t}}
 \sum\limits_{\stackrel{b_{2}=1}{b_{1}-b_{2}+n\equiv
0(p^{t-1})}}^{p^{t}}h(p^{t},b_{2}).$\\
Again, by (8.3) and (8.5),\\
$M_{A}=M_{A1}+M_{A2},$\\
where\\
 $M_{A1}=(1-p^{t-k})^{2}\Phi_{5},$\\
$M_{A2}=(1-p^{t-k})(\Phi_{6}-\Phi_{5}),$\\
and\\
$$\Phi_{6}=\#\{b_{1}, b_{2}:\ p^{t}\mid b_{1}(b_{1}+1),\
b_{1}-b_{2}+n\equiv 0(p^{t-1}),\ b_{1},b_{2}=1,2,...,p^{t} \}=2p.$$
Hence
$$M_{A2}=(1-p^{t-k})(2p-\Phi_{5}),$$
and\\
$M_{A}=(1-p^{t-k})^{2}\Phi_{5}+(1-p^{t-k})(2p-\Phi_{5})=2p(1-p^{t-k})-p^{t-k}(1-p^{t-k})\Phi_{5}.$\\
Write\\
$$
\Phi_{7}=\#\{b_{1}, b_{2}:\ p^{t}\dagger b_{1}(b_{1}+1),\
b_{1}-b_{2}+n\equiv 0(p^{t-1}),\ b_{1},b_{2}=1,2,...,p^{t} \},
$$
$$
\Phi_{8}=\#\{b_{1}, b_{2}:\ p^{t}\dagger b_{1}(b_{1}+1),\ p^{t}\mid
b_{2}(b_{2}+1),\  b_{1}-b_{2}+n\equiv 0(p^{t-1}),\
b_{1},b_{2}=1,2,...,p^{t} \},
$$
then\\
$$\Phi_{7}=p(p^{t}-2),\ \ \Phi_{8}=2p-\Phi_{5}.$$
By (8.3),\\
$M_{B}=M_{B1}+M_{B2},$\\
where\\
 $M_{B1}=(1-p^{t-k})\Phi_{8},\ \ \ M_{B2}=\Phi_{7}-\Phi_{8}.$\\
Hence\\
$M_{B}=(1-p^{t-k})\Phi_{8}+\Phi_{7}-\Phi_{8}=p(p^{t}-2)-p^{t-k}\Phi_{8}
=p(p^{t}-2)-p^{t-k}(2p-\Phi_{5}). $\\
We obtain
$$M_{B}=p(p^{t}-2)-2p^{t-k+1}+p^{t-k}\Phi_{5}, $$
and\\
$M=2p(1-p^{t-k})-p^{t-k}(1-p^{t-k})\Phi_{5}+p(p^{t}-2)-2p^{t-k+1}+p^{t-k}\Phi_{5}\\
=p^{t+1}-4p^{t-k+1}+p^{2t-2k} \Phi_{5}$,\\
the lemma follows.\\
{\bf LEMMA 8.6}\ We have
\[ \Phi_{2}=\left\{\begin{array}{cc}
1,&\ \ {\rm if}\  p^{t}\mid n(n+1),p^{t}\dagger n(n-1),\\
1,&\ \ {\rm if}\ p^{t}\dagger n(n+1),p^{t}\mid n(n-1),\\
2,&\ \ {\rm if}\ p^{t}\mid n(n+1),p^{t}\mid n(n-1),\\
0,&\  \ {\rm if}\ p^{t}\dagger n(n+1),p^{t}\dagger n(n-1).
\end{array}
\right.
\]
Proof.\ By (8.4),
$$ \Phi_{2}=\#\{b:\ p^{t}\mid b(b+1),\ p^{t}\mid
(b+n)(b+n+1),\ b=1,2,...,p^{t} \},$$ the condition $p^{t}\mid
b(b+1)$ means $b=p^{t},p^{t}-1$.\ If $p^{t}\mid n(n+1)$,then
$b=p^{t}$ satisfies the condition $p^{t}\mid (b+n)(b+n+1)$; If
$p^{t}\mid n(n-1)$,then $b=p^{t}-1$ satisfies the condition
$p^{t}\mid (b+n)(b+n+1)$; If $p^{t}\dagger n(n+1)$ and $p^{t}\dagger
n(n-1)$, then neither $b=p^{t}$ nor $b=p^{t}-1$ satisfies the
condition $p^{t}\mid (b+n)(b+n+1)$. The
lemma follows.\\
{\bf LEMMA 8.7}\ We have\\
if $t=1$, then $\Phi_{5}=4$; if $1<t\leq k$, then
\[ \Phi_{5}=\left\{\begin{array}{cc}
2,&\ \ {\rm if}\ p^{t-1}\mid n,\\
2,&\ \ {\rm if}\ p^{t-1}\dagger n,p^{t-1}\mid (n+1),p^{t-1}\mid (n-1),\\
1,&\ \ {\rm if}\ p^{t-1}\dagger n,p^{t-1}\mid (n+1),p^{t-1}\dagger (n-1),\\
1,&\ \ {\rm if}\ p^{t-1}\dagger n,p^{t-1}\dagger (n+1),p^{t-1}\mid (n-1),\\
0,&\ \ {\rm if}\ p^{t-1}\dagger n,p^{t-1}\dagger
(n+1),p^{t-1}\dagger (n-1).
\end{array}
\right.
\]
Proof.\ By (8.5),
$$ \Phi_{5}=\#\{b_{1}, b_{2}:\ p^{t}\mid
b_{1}(b_{1}+1),\ p^{t}\mid b_{2}(b_{2}+1), b_{1}-b_{2}+n\equiv
0(p^{t-1}),\ b_{1},b_{2}=1,2,...,p^{t} \},$$
 we need to count the number of $b_{1},b_{2}$ that satisfy the
conditions ($\ast$):$\ p^{t}\mid b_{1}(b_{1}+1),\ p^{t}\mid
b_{2}(b_{2}+1),
b_{1}-b_{2}+n\equiv 0(p^{t-1})$.\\
We have $b_{1},b_{2}=p^{t},p^{t}-1$.\\
If $t=1$, then $\Phi_{5}=4$ is valid. We suppose that $t>1$.\\
If $p^{t-1}\mid n$, then $b_{1}=b_{2}=p^{t}$ and
$b_{1}=b_{2}=p^{t}-1$ satisfy the conditions ($\ast$);\\
If $p^{t-1}\dagger n,\ p^{t-1}\mid (n+1),p^{t-1}\mid (n-1)$, then
$b_{1}=p^{t},b_{2}=p^{t}-1$ and $b_{1}=p^{t}-1, b_{2}=p^{t}$ satisfy
the conditions ($\ast$);\\
If $p^{t-1}\dagger n,\ p^{t-1}\mid (n+1),p^{t-1}\dagger (n-1)$ ,
then only $b_{1}=p^{t},b_{2}=p^{t}-1$  satisfies the conditions
($\ast$); if $p^{t-1}\dagger n,\ p^{t-1}\dagger (n+1),p^{t-1}\mid
(n-1)$ then
only $b_{1}=p^{t}-1,b_{2}=p^{t}$  satisfies the conditions ($\ast$);\\
if $p^{t-1}\dagger n,\ p^{t-1}\dagger (n+1),p^{t-1}\dagger (n-1)$
then there are not  $b_{1}, b_{2}$  satisfy the conditions ($\ast$).\\
The lemma follows.\\
 {\bf LEMMA 8.8}\ Suppose that $1\leq t\leq k$. Then
$$J(p^{t})=p^{3t-2k-1}(p\Phi_{2}-\Phi_{5}).$$
Proof.\ By (8.2), Lemma 8.3 and Lemma 8.4,\\
$J(p)=J_{1}(p)-J_{2}(p)=p^{2}-4p^{2-k}+p^{3-2k}\Phi_{2}-p^{2}(1-2p^{-k})^{2}
=p^{2-2k}(p\Phi_{2}-4).$\\
We suppose $1<t\leq k$.  By (8.2), Lemma 8.3 and Lemma 8.5,\\
$
J(p^{t})=J_{1}(p^{t})-J_{2}(p^{t})=p^{2t}-4p^{2t-k}+p^{3t-2k}\Phi_{2}-p^{t-1}(p^{t+1}-4p^{t-k+1}+p^{2t-2k}
\Phi_{5})\\
= p^{3t-2k-1}(p\Phi_{2}-\Phi_{5}),$\\
the lemma follows.\\
By Lemma 8.1, Lemma 8.2 and Lemma 8.8,\\
$\mathfrak{S}(n)=\varrho^{2}\prod\limits_{p}(1+\sum\limits_{t=1}^{k}p^{-2t}f(p^{t})^{2}J(p^{t}))
=\varrho^{2}\prod\limits_{p}(1+\sum\limits_{t=1}^{k}p^{-2t}(\frac{p^{k}}{p^{k}-2})^{2}
p^{3t-2k-1}(p\Phi_{2}-\Phi_{5}))\\
=\varrho^{2}\prod\limits_{p}(1+(p^{k}-2)^{-2}\sum\limits_{t=1}^{k}p^{t-1}(p\Phi_{2}-\Phi_{5})),$\\
that is\\
\bi \mathfrak{S}(n)=\varrho^{2}\prod\limits_{p}\Delta_{p},\ei where
\bi
\Delta_{p}=1+(p^{k}-2)^{-2}\sum\limits_{t=1}^{k}p^{t-1}(p\Phi_{2}-\Phi_{5}).
\ei
 Write\\
\bi \Pi_{1}=\prod\limits_{p\mid n(n+1)}\Delta_{p},\ \
\Pi_{1A}=\prod\limits_{\stackrel{p\mid n(n+1)}{p\mid
n(n-1)}}\Delta_{p},\ \ \Pi_{1B}=\prod\limits_{\stackrel{p\mid
n(n+1)}{p\dagger n(n-1)}}\Delta_{p}, \ \ \Pi_{1}=\Pi_{1A}\Pi_{1B},
\ei
\bi
 \Pi_{2}=\prod\limits_{p\dagger n(n+1)}\Delta_{p},\ \
\Pi_{2A}=\prod\limits_{\stackrel{p\dagger n(n+1)}{p\mid
n-1}}\Delta_{p},\ \ \Pi_{2B}=\prod\limits_{\stackrel{p\dagger
n(n+1)}{p\dagger n-1}}\Delta_{p},\ \ \ \
\Pi_{2}=\Pi_{2A}\Pi_{2B},\ei
 then
 \bi \mathfrak{S}(n)=\varrho^{2}\Pi_{1A}\Pi_{1B}\Pi_{2A}\Pi_{2B}. \ei
For fixed $n,p$, we write\\
$$p^{\beta_{1}}\parallel n,\ p^{\beta_{2}}\parallel n+1,\ p^{\beta_{3}}\parallel n-1.$$
 We consider $\Pi_{1A}$ first. \\
{\bf LEMMA 8.9}\ We have\\
$\Pi_{1A} =\prod\limits_{\stackrel{p\mid n}{p^{k}\dagger
n}}(1-4(p^{k}-2)^{-2})\prod\limits_{p^{k}\mid n}(1+2(p^{k}-2)^{-1})
 \prod\limits_{\stackrel{2^{k}\mid
n+1}{or 2^{k}\mid
n-1}}(1+(2^{k}-4)(2^{k}-2)^{-2})\times\\
\prod\limits_{\stackrel{2^{k}\dagger
n+1}{2^{k}\dagger n-1, 2\dagger n}}(1-4(2^{k}-2)^{-2}).$\\
Proof.\ We have\\
\bi
\Pi_{1A}=\prod\limits_{\stackrel{p^{k}\dagger n}{p\mid
n}}\Delta_{p}\prod\limits_{p^{k}\mid
n}\Delta_{p}\prod\limits_{2\dagger n}\Delta_{2}.
\ei
If $p\mid
n(n+1),\ p\mid n(n-1)$ and $ p\dagger n $,  then
$p=2,\beta_{1}=0,\ \min(\beta_{2}, \beta_{3})=1$.\\
Hence, by Lemma 8.6 and Lemma 8.7, \\
 $\Phi_{2}(2)=2, \Phi_{5}(2)=4$, and then
 $2\Phi_{2}(2)-\Phi_{5}(2)=0.$\\
When $1<t\leq k, 2\dagger n,$
\[ \Phi_{2}(2^{t})=\left\{\begin{array}{cc}
1,&\ \ {\rm if}\  t\leq \max(\beta_{2},\beta_{3}),\\
0,&\  \ {\rm if}\ t>\max(\beta_{2},\beta_{3}).
\end{array}
\right.
\]
 \[ \Phi_{5}(2^{t})=\left\{\begin{array}{cc}
2,&\ \ {\rm if}\ t=2,\\
1,&\ \ {\rm if}\ 2<t\leq \max(\beta_{2},\beta_{3})+1,\\
0,&\ \ {\rm if}\ t>\max(\beta_{2},\beta_{3})+1.
\end{array}
\right.
\]
Write $\beta=\max(\beta_{2},\beta_{3})\geq 2$.   We have \\
\[ 2\Phi_{2}(2^{t})-\Phi_{5}(2^{t})=\left\{\begin{array}{cc}
0,&\  \ {\rm if}\ t=2,\\
1,&\ \ {\rm if}\  2<t\leq \beta,\\
-1,&\ \ {\rm if}\  t=1+\beta,\\
0,&\  \ {\rm if}\ 1+\beta<t\leq k.
\end{array}
\right.
\]
Hence
\[ \Delta_{2}=\left\{\begin{array}{cc}
1+(2^{k}-2)^{-2}(\sum\limits_{3\leq t\leq \beta}2^{t-1}-2^{\beta}),&\  \ {\rm if}\ \beta<k,\\
1+(2^{k}-2)^{-2}\sum\limits_{3\leq t\leq k}2^{t-1},&\  \ {\rm if}\
\beta\geq k.
\end{array}
\right.
\]
Then
\bi
\Delta_{2}=\left\{\begin{array}{cc}
1-4(2^{k}-2)^{-2},&\  \ {\rm if}\ \beta<k,\\
1+(2^{k}-4)(2^{k}-2)^{-2},&\  \ {\rm if}\ \beta\geq k.
\end{array}
\right.
 \ei
 We suppose that $p\mid n(n+1),\ p\mid n(n-1)$ and $ p\mid n $, then
$\beta_{2}=\beta_{3}=0$.\\
By Lemma 8.6 and Lemma 8.7,\\
$$p\Phi_{2}(p)-\Phi_{5}(p)=2p-4,$$
\[ \Phi_{2}(p^{t})=\left\{\begin{array}{cc}
2,&\  \ {\rm if}\ p^{t}\mid n,\\
0,&\  \ {\rm if}\ p^{t}\dagger n,
\end{array}
\right.
\]
and for $1<t\leq k$,
\[ \Phi_{5}(p^{t})=\left\{\begin{array}{cc}
2,&\ \ {\rm if}\ p^{t-1}\mid n,\\
0,&\ \ {\rm if}\ p^{t-1}\dagger n.
\end{array}
\right.
\]
Hence
\[ p\Phi_{2}(p^{t})-\Phi_{5}(p^{t})=\left\{\begin{array}{cc}
2p-2,&\  \ {\rm if}\ t\leq \beta_{1}\\
-2,&\ \ {\rm if}\  t=1+\beta_{1},\\
0,&\  \ {\rm if}\ 1+\beta_{1}<t\leq k.
\end{array}
\right.
\]
Then
\[ \Delta_{p}=\left\{\begin{array}{cc}
1+(p^{k}-2)^{-2}(2p-4)+(p^{k}-2)^{-2}(\sum\limits_{2\leq t\leq \beta_{1}}p^{t-1}(2p-2)-2p^{\beta_{1}}),&\  \ {\rm if}\ \beta_{1}<k,\\
1+(p^{k}-2)^{-2}(2p-4)+(p^{k}-2)^{-2}(\sum\limits_{2\leq t\leq
k}p^{t-1}(2p-2)),&\  \ {\rm if}\ \beta_{1}\geq k.
\end{array}
\right.
\]
Consequently
\bi
\Delta_{p}=\left\{\begin{array}{cc}
1-4(p^{k}-2)^{-2},&\  \ {\rm if}\ 1\leq\beta_{1}<k,\\
1+2(p^{k}-2)^{-1},&\ \ {\rm if}\ \beta_{1}\geq k.
\end{array}
\right.
 \ei
By (8.11), (8.12) and (8.13), the lemma follows.\\
Secondly, we consider $\Pi_{1B}$.\\
{\bf LEMMA 8.10}\ We have\\
$\Pi_{1B} =\prod\limits_{\stackrel{p^{k}\dagger n+1}{p\mid n+1,
p\dagger n-1}}(1-4(p^{k}-2)^{-2})\prod\limits_{\stackrel{p^{k}\mid
n+1}{p\dagger n-1}}(1+(p^{k}-4)(p^{k}-2)^{-2}).$\\
Proof.\ We have  $ p\mid n(n+1), p\dagger n(n-1),
\beta_{1}=\beta_{3}=0, \beta_{2} \geq 1$, and
 by Lemma 8.6 and Lemma 8.7, \\
$$p\Phi_{2}(p)-\Phi_{5}(p)=p-4, $$
for $1<t\leq k$,\\
\[ \Phi_{2}(p^{t})=\left\{\begin{array}{cc}
1,&\ \ {\rm if}\  p^{t}\mid n+1,\\
0,&\  \ {\rm if}\ p^{t}\dagger n+1,
\end{array}
\right.
\]
\[ \Phi_{5}(p^{t})=\left\{\begin{array}{cc}
1,&\ \ {\rm if}\  p^{t-1}\mid n+1,\\
0,&\  \ {\rm if}\ p^{t-1}\dagger n+1,
\end{array}
\right.
\]
then
\[ p\Phi_{2}(p^{t})-\Phi_{5}(p^{t})=\left\{\begin{array}{cc}
p-1,&\  \ {\rm if}\ t\leq \beta_{2}\\
-1,&\ \ {\rm if}\  t=1+\beta_{2},\\
0,&\  \ {\rm if}\ 1+\beta_{2}<t\leq k,
\end{array}
\right.
\]
and
\[ \Delta_{p}=\left\{\begin{array}{cc}
1+(p^{k}-2)^{-2}(p-4)+(p^{k}-2)^{-2}(\sum\limits_{2\leq t\leq \beta_{2}}p^{t-1}(p-1)-p^{\beta_{2}}),&\  \ {\rm if}\ \beta_{2}<k,\\
1+(p^{k}-2)^{-2}(p-4)+(p^{k}-2)^{-2}(\sum\limits_{2\leq t\leq
k}p^{t-1}(p-1)),&\  \ {\rm if}\ \beta_{2}\geq k.
\end{array}
\right.
\]
Consequently
\bi
\Delta_{p}=\left\{\begin{array}{cc}
1-4(p^{k}-2)^{-2},&\  \ {\rm if}\ 1\leq\beta_{2}<k,\\
1+(p^{k}-4)(p^{k}-2)^{-2},&\ \ {\rm if}\ \beta_{2}\geq k.
\end{array}
\right.
 \ei
By (8.8) and (8.14), the lemma follows.\\
Now, we consider $\Pi_{2}$.  By $ p\dagger n(n+1)$, Lemma 8.6 and Lemma 8.7, \\
\bi
 \Phi_{2}(p^{t})=\left\{\begin{array}{cc}
1,&\ \ {\rm if}\  p^{t}\mid n-1,\\
0,&\  \ {\rm if}\ p^{t}\dagger n-1,
\end{array}
\right.
\ei
and for $1<t\leq k$
 \bi
 \Phi_{5}(p^{t})=\left\{\begin{array}{cc}
1,&\ \ {\rm if}\  p^{t-1}\mid n-1,\\
0,&\  \ {\rm if}\ p^{t-1}\dagger n-1.
\end{array}
\right.
\ei
 {\bf LEMMA 8.11}\ We have\\
$\Pi_{2A} =\prod\limits_{\stackrel{p\dagger n(n+1)}{p^{k}\dagger
n-1, p\mid n-1}}(1-4(p^{k}-2)^{-2})\prod\limits_{\stackrel{p\dagger
n(n+1)}{p^{k}\mid n-1}}(1+(p^{k}-4)(p^{k}-2)^{-2}).$\\
Proof.\ We have  $ p\dagger n(n+1), p\mid (n-1),
\beta_{1}=\beta_{2}=0, \beta_{3} \geq 1$, and
 by Lemma 8.6, Lemma 8.7, (8.15) and (8.16), \\
$ p\Phi_{2}(p)-\Phi_{5}(p)=p-4,$\\
and for $1<t\leq k$, by (8.15) and (8.16),\\
\[ p\Phi_{2}(p^{t})-\Phi_{5}(p^{t})=\left\{\begin{array}{cc}
p-1,&\  \ {\rm if}\ t\leq \beta_{3}\\
-1,&\ \ {\rm if}\  t=1+\beta_{3},\\
0,&\  \ {\rm if}\ 1+\beta_{3}<t\leq k,
\end{array}
\right.
\]
and
\[ \Delta_{p}=\left\{\begin{array}{cc}
1+(p^{k}-2)^{-2}(p-4)+(p^{k}-2)^{-2}(\sum\limits_{2\leq t\leq \beta_{3}}p^{t-1}(p-1)-p^{\beta_{3}}),&\  \ {\rm if}\ \beta_{3}<k,\\
1+(p^{k}-2)^{-2}(p-4)+(p^{k}-2)^{-2}(\sum\limits_{2\leq t\leq
k}p^{t-1}(p-1)),&\  \ {\rm if}\ \beta_{3}\geq k.
\end{array}
\right.
\]
Consequently
\bi
\Delta_{p}=\left\{\begin{array}{cc}
1-4(p^{k}-2)^{-2},&\  \ {\rm if}\ 1\leq\beta_{3}<k,\\
1+(p^{k}-4)(p^{k}-2)^{-2},&\ \ {\rm if}\ \beta_{3}\geq k.
\end{array}
\right.
 \ei
By (8.9) and (8.17), the lemma follows.\\
{\bf LEMMA 8.12}\ We have\\
$\Pi_{2B} =\prod\limits_{p\dagger
n(n+1)(n-1)}(1-4(p^{k}-2)^{-2}).$\\
Proof.\ We have  $ p\dagger n(n+1)(n-1),
\beta_{1}=\beta_{2}=\beta_{3}=0$, and
 by Lemma 8.6, Lemma 8.7, (8.15) and (8.16), \\
$ p\Phi_{2}(p)-\Phi_{5}(p)=-4,$\\
and for $1<t\leq k$,\\
$ \Phi_{2}(p^{t})=0,\  \Phi_{5}(p^{t})=0,$\\
and\\
$ p\Phi_{2}(p^{t})-\Phi_{5}(p^{t})=0$, and
 $\Delta_{p}=1-4(p^{k}-2)^{-2}$,\\
 the lemma follows.\\
{\bf LEMMA 8.13}\ We have\\
$\mathfrak{S}(n) =\varrho^{2}\prod\limits_{\stackrel{p\mid
n}{p^{k}\dagger n}}(1-4(p^{k}-2)^{-2})\prod\limits_{p^{k}\mid
n}(1+2(p^{k}-2)^{-1})
 \prod\limits_{\stackrel{2^{k}\mid
n+1}{or 2^{k}\mid
n-1}}(1+(2^{k}-4)(2^{k}-2)^{-2})\times\\
\prod\limits_{\stackrel{2^{k}\dagger n+1}{2^{k}\dagger n-1, 2\dagger
n}}(1-4(2^{k}-2)^{-2}) \prod\limits_{\stackrel{p^{k}\dagger
n+1}{p\mid n+1, p\dagger
n-1}}(1-4(p^{k}-2)^{-2})\prod\limits_{\stackrel{p^{k}\mid
n+1}{p\dagger n-1}}(1+(p^{k}-4)(p^{k}-2)^{-2})\times\\
\prod\limits_{\stackrel{p\dagger n(n+1)}{p^{k}\dagger n-1, p\mid
n-1}}(1-4(p^{k}-2)^{-2})\prod\limits_{\stackrel{p\dagger
n(n+1)}{p^{k}\mid n-1}}(1+(p^{k}-4)(p^{k}-2)^{-2})
\prod\limits_{p\dagger n(n+1)(n-1)}(1-4(p^{k}-2)^{-2}).$\\
Proof.\ By Lemma 8.9, Lemma 8.10, Lemma 8.11, Lemma 8.12 and (8.10),
the lemma follows.\\
{\bf LEMMA 8.14}\ We have\\
$\mathfrak{S}(n)
=\varrho^{2}\prod\limits_{p\not=2}(1-4(p^{k}-2)^{-2})
\prod\limits_{\stackrel{2^{k}\mid n+1}{or 2^{k}\mid
n-1}}(1+(2^{k}-4)(2^{k}-2)^{-2})\prod\limits_{\stackrel{2^{k}\dagger
n+1}{2^{k}\dagger n-1, 2\dagger n}}(1-4(2^{k}-2)^{-2})\times\\
\prod\limits_{\stackrel{2\mid n}{2^{k}\dagger
n}}(1-4(2^{k}-2)^{-2})\prod\limits_{2^{k}\mid
n}(1+2(2^{k}-2)^{-1})\prod\limits_{\stackrel{p^{k}\mid
n}{p\not=2}}\frac{p^{k}-2}{p^{k}-4}
\prod\limits_{\stackrel{p^{k}\mid
n+1}{p\not=2}}\frac{p^{k}-3}{p^{k}-4}\prod\limits_{\stackrel{p^{k}\mid
n-1}{p\not=2}}\frac{p^{k}-3}{p^{k}-4}.$\\
Proof.\ By Lemma 8.13 and \\
$\prod\limits_{p\dagger
n(n+1)(n-1)}(1-4(p^{k}-2)^{-2})=\prod\limits_{p\not=2}(1-4(p^{k}-2)^{-2})
\prod\limits_{\stackrel{p\mid
n(n+1)(n-1)}{p\not=2}}(1-4(p^{k}-2)^{-2})^{-1},$\\
the lemma follows at once.\\
\section{ The function $G(u;s)$ }
\setcounter{equation}{0}
By (7.4),
 \bi
 S_{E}(x,Q)=\frac{1}{2\pi i}\int_{(2)}\sum\limits_{u\leq Q}\sum\limits_{v=1}^{\infty}\mathfrak{S}(uv)
(uv)^{-s}\frac{x^{s+1}}{s(s+1)}ds,
 \ei
where $\int_{(2)}$ means $\int_{2-i\infty}^{2+i\infty}$.\\
 Write
$$
G(u;s)=\sum\limits_{v=1}^{\infty}\mathfrak{S}(uv)v^{-s},\ \ \  \Re
s>1,
 $$
then
\bi
 S_{E}(x,Q)=\sum\limits_{u\leq Q}\frac{1}{2\pi
i}\int_{(2)}u^{-s}G(u;s)\frac{x^{s+1}}{s(s+1)}ds.
 \ei
Let
$$
G_{1}(u,s)=\sum\limits_{\stackrel{v=1}{2^{k}\mid
uv}}^{\infty}\mathfrak{S}(uv)v^{-s},
 \ G_{2}(u,s)=\sum\limits_{\stackrel{v=1}{2\mid uv, 2^{k}\dagger
 uv}}^{\infty}\mathfrak{S}(uv)v^{-s},$$
$$ G_{3}(u,s)=\sum\limits_{\stackrel{2\dagger uv,v=1}{2^{k}\dagger
uv+1,2^{k}\dagger uv-1 }}^{\infty}\mathfrak{S}(uv)v^{-s}, \
G_{4}(u,s)=\sum\limits_{\stackrel{2\dagger uv,v=1}{2^{k}\mid uv+1 or
2^{k}\mid uv-1 }}^{\infty}\mathfrak{S}(uv)v^{-s},
 $$
 then
 \bi
 G(u;s)= G_{1}(u,s)+G_{2}(u,s)+G_{3}(u,s)+G_{4}(u,s).
 \ei
As in [3], for integer $a$, let\\
$\xi_{k}=e^{2\pi i/k }=e(\frac{1}{k}),\
\zeta(s;\xi_{k}^{a})=\sum\limits_{m=1}^{\infty}\xi_{k}^{am}m^{-s}, \
\zeta(s;a,k)=\sum\limits_{n\equiv a(k)}n^{-s}, \Re s>1$,\\
and for $ 0<a\leq 1$, let $
\zeta(s,a)=\sum\limits_{n=0}^{\infty}(n+a)^{-s}, \Re s>1. $\\
When $1\leq a\leq k$, we have
$\zeta(s;a,k)=k^{-s}\zeta(s,\frac{a}{k})$.\\
For odd integer, we define the multiplicative functions $g_{1}(d),
g_{2}(d)$ as follows
\bi
 g_{1}(p^{t})=\frac{2}{p^{k}-4}, \ \
g_{2}(p^{t})=\frac{1}{p^{k}-4}, \ t\geq 1, p>2,
\ei
then
\bi \mu(d)^{2}g_{1}(d)\ll \tau(d)d^{-k}, \ \  \
 \mu(d)^{2}g_{2}(d)\ll d^{-k}. \ei
 For fixed $u,\ d_{1}, d_{2},\ d_{3}$ and $\Re s>1$, \ write\\
$H_{1}(s)=H_{1}(s;u,d_{1},d_{2},d_{3})=\sum\limits_{\stackrel{2^{k}\mid
uv, d_{1}^{k}\mid uv}{d_{2}^{k}\mid uv+1, d_{3}^{k}\mid
uv-1}}v^{-s}, $\\
$H_{2}(s)=H_{2}(s;u,d_{1},d_{2},d_{3})=\sum\limits_{\stackrel{2\mid
uv, d_{1}^{k}\mid uv}{d_{2}^{k}\mid uv+1, d_{3}^{k}\mid
uv-1}}v^{-s},$\\
$H_{3}(s)=H_{3}(s;u,d_{1},d_{2},d_{3})=\sum\limits_{\stackrel{
d_{1}^{k}\mid uv}{d_{2}^{k}\mid uv+1, d_{3}^{k}\mid
uv-1}}v^{-s},$\\
$H_{4}(s)=H_{4}(s;u,d_{1},d_{2},d_{3})=\sum\limits_{\stackrel{2^{k}\mid
uv+1, d_{1}^{k}\mid uv}{d_{2}^{k}\mid uv+1, d_{3}^{k}\mid
uv-1}}v^{-s},$\\
$H_{5}(s)=H_{5}(s;u,d_{1},d_{2},d_{3})=\sum\limits_{\stackrel{2^{k}\mid
uv-1, d_{1}^{k}\mid uv}{d_{2}^{k}\mid uv+1, d_{3}^{k}\mid
uv-1}}v^{-s}$.\\
By Lemma 8.14 and (9.4),\\
$\mathfrak{S}(n)
=\varrho^{2}\prod\limits_{p\not=2}(1-4(p^{k}-2)^{-2})
\prod\limits_{\stackrel{2^{k}\mid n+1}{or 2^{k}\mid
n-1}}(1+(2^{k}-4)(2^{k}-2)^{-2})\times\\
\prod\limits_{\stackrel{2^{k}\dagger n+1}{2^{k}\dagger n-1, 2\dagger
n}}(1-4(2^{k}-2)^{-2}) \prod\limits_{\stackrel{2\mid n}{2^{k}\dagger
n}}(1-4(2^{k}-2)^{-2})\prod\limits_{2^{k}\mid
n}(1+2(2^{k}-2)^{-1})\times\\
\sum\limits_{\stackrel{d_{1}^{k}\mid
n}{(d_{1},2)=1}}\mu(d_{1})^{2}g_{1}(d_{1})
\sum\limits_{\stackrel{d_{2}^{k}\mid
n+1}{(d_{2},2)=1}}\mu(d_{2})^{2}g_{2}(d_{2})\sum\limits_{\stackrel{d_{3}^{k}\mid
n-1}{(d_{3},2)=1}}\mu(d_{3})^{2}g_{2}(d_{3}).$\\
If we write \\
$$h(d_{1},d_{2},d_{3})=\mu(d_{1}d_{2}d_{3})^{2}g_{1}(d_{1})g_{2}(d_{2})g_{2}(d_{3}),$$
then, by (9.5),
 \bi
  h(d_{1},d_{2},d_{3})\ll
\tau(d_{1})d_{1}^{-k}d_{2}^{-k}d_{3}^{-k},
\ei
 and
\bi
\left.
\begin{array}{cc}
\mathfrak{S}(n)=\varrho^{2}\prod\limits_{p\not=2}(1-4(p^{k}-2)^{-2})
\prod\limits_{\stackrel{2^{k}\mid n+1}{or 2^{k}\mid
n-1}}(1+(2^{k}-4)(2^{k}-2)^{-2})\times\\
\prod\limits_{\stackrel{2^{k}\dagger n+1}{2^{k}\dagger n-1, 2\dagger
n}}(1-4(2^{k}-2)^{-2}) \prod\limits_{\stackrel{2\mid n}{2^{k}\dagger
n}}(1-4(2^{k}-2)^{-2})\prod\limits_{2^{k}\mid
n}(1+2(2^{k}-2)^{-1})\times\\
\sum\limits_{\stackrel{(d_{1}d_{2}d_{3},2)=1}{d_{1}^{k}\mid
n,d_{2}^{k}\mid n+1,d_{3}^{k}\mid n-1}}h(d_{1},d_{2},d_{3}).
\end{array}
\right.
\ei
 {\bf LEMMA 9.1}\ For $\Re s>1$, we have\\
$G_{1}(u,s)=\varrho^{2}\prod\limits_{p\not=2}(1-4(p^{k}-2)^{-2})(1+2(2^{k}-2)^{-1})\sum\limits_{(d_{1}d_{2}d_{3},2)=1}
h(d_{1},d_{2},d_{3})H_{1}(s),\\
G_{2}(u,s)=\varrho^{2}\prod\limits_{p\not=2}(1-4(p^{k}-2)^{-2})(1-4(2^{k}-2)^{-2})\sum\limits_{(d_{1}d_{2}d_{3},2)=1}
h(d_{1},d_{2},d_{3})(H_{2}(s)-H_{1}(s)),\\
G_{3}(u,s)=\varrho^{2}\prod\limits_{p\not=2}(1-4(p^{k}-2)^{-2})(1-4(2^{k}-2)^{-2})\sum\limits_{(d_{1}d_{2}d_{3},2)=1}
h(d_{1},d_{2},d_{3})(H_{3}(s)-H_{2}(s)-H_{4}(s)-H_{5}(s)),\\
G_{4}(u,s)=\varrho^{2}\prod\limits_{p\not=2}(1-4(p^{k}-2)^{-2})(1+(2^{k}-4)(2^{k}-2)^{-2})\sum\limits_{(d_{1}d_{2}d_{3},2)=1}
h(d_{1},d_{2},d_{3})(H_{4}(s)+H_{5}(s)),$\\
$G(u;s)=\varrho^{2}\prod\limits_{p\not=2}(1-4(p^{k}-2)^{-2})
\sum\limits_{(d_{1}d_{2}d_{3},2)=1}
h(d_{1},d_{2},d_{3})(2^{k+1}(2^{k}-2)^{-2}H_{1}(s)+\\
(1-4(2^{k}-2)^{-2})H_{3}(s)+2^{k}(2^{k}-2)^{-2}H_{4}(s)+
2^{k}(2^{k}-2)^{-2}H_{5}(s)).$\\
Proof.\ By (9.3) and (9.7), the lemma follows.\\
We only consider the case: $\mu(d_{1}d_{2}d_{3})^{2}=1,
\ (d_{2}d_{3},u)=1 $. \ Write\\
\bi \frac{2^{k}d_{1}^{k}u}{(u,2^{k}d_{1}^{k})}b_{1}\equiv 1 (mod\
d_{2}^{k}), \ \frac{2^{k}d_{1}^{k}u}{(u,2^{k}d_{1}^{k})}b_{2}\equiv
1 (mod\ d_{3}^{k}), \ei
 \bi \left\{\begin{array}{cc}
d_{3}^{k}\overline{d_{3}^{k}}\equiv 1 (mod\ d_{2}^{k}),\\
d_{2}^{k}\overline{d_{2}^{k}}\equiv 1 (mod\ d_{3}^{k}),
\end{array}
\right.
 \ei
\bi A_{1}\equiv
-b_{1}d_{3}^{k}\overline{d_{3}^{k}}+b_{2}d_{2}^{k}\overline{d_{2}^{k}}
\ (mod\ d_{2}^{k}d_{3}^{k}), \ 1\leq A_{1}\leq d_{2}^{k}d_{3}^{k},
(A_{1},d_{2}d_{3})=1, \ei
  \bi \frac{d_{1}^{k}u}{(u,d_{1}^{k})}b_{5}\equiv 1 (mod\ d_{2}^{k}),
\ \frac{d_{1}^{k}u}{(u,d_{1}^{k})}b_{6}\equiv 1 (mod\ d_{3}^{k}),
\ei
\bi A_{3}\equiv
-b_{5}d_{3}^{k}\overline{d_{3}^{k}}+b_{6}d_{2}^{k}\overline{d_{2}^{k}}
\ (mod\ d_{2}^{k}d_{3}^{k}), \ 1\leq A_{3}\leq d_{2}^{k}d_{3}^{k},
(A_{3},d_{2}d_{3})=1, \ei
 \bi \frac{d_{1}^{k}u}{(u,d_{1}^{k})}b_{7}\equiv 1 (mod\
2^{k}d_{2}^{k}), \ \frac{d_{1}^{k}u}{(u,d_{1}^{k})}b_{8}\equiv 1
(mod\ d_{3}^{k}), \ei
 \bi \left\{\begin{array}{cc}
d_{3}^{k}\overline{(d_{3}^{k})}\equiv 1 (mod\ 2^{k}d_{2}^{k}),\\
2^{k}d_{2}^{k}\overline{(2^{k}d_{2}^{k})}\equiv 1 (mod\ d_{3}^{k}),
\end{array}
\right.
 \ei
\bi A_{4}\equiv
-b_{7}d_{3}^{k}\overline{(d_{3}^{k})}+b_{8}2^{k}d_{2}^{k}\overline{(2^{k}d_{2}^{k})}
\ (mod\ 2^{k}d_{2}^{k}d_{3}^{k}), \ 1\leq A_{4}\leq
2^{k}d_{2}^{k}d_{3}^{k}, (A_{4},2d_{2}d_{3})=1, \ei \bi
\frac{d_{1}^{k}u}{(u,d_{1}^{k})}b_{9}\equiv 1 (mod\ d_{2}^{k}), \
\frac{d_{1}^{k}u}{(u,d_{1}^{k})}b_{10}\equiv 1 (mod\
2^{k}d_{3}^{k}), \ei \bi \left\{\begin{array}{cc}
2^{k}d_{3}^{k}\overline{\{2^{k}d_{3}^{k}\}}\equiv 1 (mod\ d_{2}^{k}),\\
d_{2}^{k}\overline{\{d_{2}^{k}\}}\equiv 1 (mod\ 2^{k}d_{3}^{k}),
\end{array}
\right.
 \ei
\bi A_{5}\equiv
-b_{9}2^{k}d_{3}^{k}\overline{\{2^{k}d_{3}^{k}\}}+b_{10}d_{2}^{k}\overline{\{d_{2}^{k}\}}
\ (mod\ 2^{k}d_{2}^{k}d_{3}^{k}), \ 1\leq A_{5}\leq
2^{k}d_{2}^{k}d_{3}^{k}, (A_{5},2d_{2}d_{3})=1. \ei
 {\bf LEMMA 9.2}\ For $\Re s>1$, we have\\
$H_{1}(s)=\frac{(u,2^{k}d_{1}^{k})^{s}}{2^{ks}d_{1}^{ks}}\zeta(s;A_{1},d_{2}^{k}d_{3}^{k}), \ H_{3}(s)=\frac{(u,d_{1}^{k})^{s}}{d_{1}^{ks}}\zeta(s;A_{3},d_{2}^{k}d_{3}^{k}),$\\
$H_{4}(s)=\frac{(u,d_{1}^{k})^{s}}{d_{1}^{ks}}\zeta(s;A_{4},2^{k}d_{2}^{k}d_{3}^{k}),\ H_{5}(s)=\frac{(u,d_{1}^{k})^{s}}{d_{1}^{ks}}\zeta(s;A_{5},2^{k}d_{2}^{k}d_{3}^{k}).$\\
Proof.\ We have\\
$2^{k}\mid uv, d_{1}^{k}\mid uv, d_{2}^{k}\mid uv+1, d_{3}^{k}\mid
uv-1. $\\
Write  $a=\frac{2^{k}d_{1}^{k}}{(u,2^{k}d_{1}^{k})},  v=av_{1},$ then \\
$$uav_{1}\equiv -1 (mod\ d_{2}^{k}),\ uav_{1}\equiv 1 (mod\ d_{3}^{k}).$$
By (9.8) and (9.9),
\[ \left\{\begin{array}{cc}
v_{1}\equiv -b_{1}(mod\ d_{2}^{k}),\\
v_{1}\equiv b_{2}(mod\ d_{3}^{k}),
\end{array}
\right.
\]
hence, by (9.10), we obtain  $v_{1}\equiv A_{1} (mod\
d_{2}^{k}d_{3}^{k})$ and \\
$H_{1}(s)=a^{-s}\sum\limits_{v_{1}\equiv A_{1}(d_{2}^{k}d_{3}^{k})}v_{1}^{-s}=a^{-s}\zeta(s;A_{1},d_{2}^{k}d_{3}^{k}). $\\
Similarly, the formulas for $H_{3}(s), H_{4}(s) $ and $H_{5}(s)$
follow from (9.11),(9.12),(9.13),(9.14), (9.15),(9.16),(9.17) and (9.18), and the lemma follows.\\
We use the following results of [3].\\
Let $G(s)=-i(2\pi)^{s-1}\Gamma(1-s)$,  then
 \bi
\zeta(s;a,k)=G(s)k^{-s}(e^{\frac{1}{2}i\pi
s}\zeta(1-s;\xi_{k}^{a})-e^{\frac{-1}{2}i\pi
s}\zeta(1-s;\xi_{k}^{-a})). \ei
 By Lemma 9.1 and Lemma 9.2($\Re s>1$),
\bi \left.
\begin{array}{cc}
G(u;s)=\varrho^{2}\prod\limits_{p\not=2}(1-4(p^{k}-2)^{-2})
\sum\limits_{\stackrel{(d_{2}d_{3},u)=1}{(d_{1}d_{2}d_{3},2)=1}}
h(d_{1},d_{2},d_{3})\times\\
\{
2^{k+1}(2^{k}-2)^{-2}\frac{(u,2^{k}d_{1}^{k})^{s}}{2^{ks}d_{1}^{ks}}
\zeta(s;A_{1},d_{2}^{k}d_{3}^{k})\\
+
(1-4(2^{k}-2)^{-2})\frac{(u,d_{1}^{k})^{s}}{d_{1}^{ks}}\zeta(s;A_{3},d_{2}^{k}d_{3}^{k})\\
+
2^{k}(2^{k}-2)^{-2}\frac{(u,d_{1}^{k})^{s}}{d_{1}^{ks}}\zeta(s;A_{4},2^{k}d_{2}^{k}d_{3}^{k})\\
+2^{k}(2^{k}-2)^{-2}\frac{(u,d_{1}^{k})^{s}}{d_{1}^{ks}}\zeta(s;A_{5},2^{k}d_{2}^{k}d_{3}^{k})
\}.
\end{array}
\right. \ei We take (9.20) as the  analytic continuation of $G(u;s),
\Re s>\frac{1}{k}-1$,
and $G(u;s)$ has at most one simple pole at $s=1$.\\
Write\\
$
G_{(1)}(u,s)=\sum\limits_{\stackrel{(d_{2}d_{3},u)=1}{(d_{1}d_{2}d_{3},2)=1}}
h(d_{1},d_{2},d_{3})\frac{(u,2^{k}d_{1}^{k})^{s}}{2^{ks}d_{1}^{ks}}
\zeta(s;A_{1},d_{2}^{k}d_{3}^{k}),\\
G_{(3)}(u,s)=\sum\limits_{\stackrel{(d_{2}d_{3},u)=1}{(d_{1}d_{2}d_{3},2)=1}}
h(d_{1},d_{2},d_{3})\frac{(u,d_{1}^{k})^{s}}{d_{1}^{ks}}\zeta(s;A_{3},d_{2}^{k}d_{3}^{k}),\\
G_{(4)}(u,s)=\sum\limits_{\stackrel{(d_{2}d_{3},u)=1}{(d_{1}d_{2}d_{3},2)=1}}
h(d_{1},d_{2},d_{3})\frac{(u,d_{1}^{k})^{s}}{d_{1}^{ks}}\zeta(s;A_{4},2^{k}d_{2}^{k}d_{3}^{k}),\\
G_{(5)}(u,s)=\sum\limits_{\stackrel{(d_{2}d_{3},u)=1}{(d_{1}d_{2}d_{3},2)=1}}
h(d_{1},d_{2},d_{3})\frac{(u,d_{1}^{k})^{s}}{d_{1}^{ks}}\zeta(s;A_{5},2^{k}d_{2}^{k}d_{3}^{k}),
$\\
 then
 \bi \left.
\begin{array}{cc}
G(u;s)=\varrho^{2}\prod\limits_{p\not=2}(1-4(p^{k}-2)^{-2})
2^{k+1}(2^{k}-2)^{-2}G_{(1)}(u,s)\\
+\varrho^{2}\prod\limits_{p\not=2}(1-4(p^{k}-2)^{-2})
(1-4(2^{k}-2)^{-2})G_{(3)}(u,s)\\
+\varrho^{2}\prod\limits_{p\not=2}(1-4(p^{k}-2)^{-2})2^{k}(2^{k}-2)^{-2}G_{(4)}(u,s)\\
+\varrho^{2}\prod\limits_{p\not=2}(1-4(p^{k}-2)^{-2})2^{k}(2^{k}-2)^{-2}G_{(5)}(u,s).
\end{array}
\right. \ei
\section{ The formula for $S_{E}(x,Q)$ }
\setcounter{equation}{0}
Let $T>1$, by (9.2),\\
$
 S_{E}(x,Q)=\sum\limits_{u\leq Q}\frac{1}{2\pi
i}(\int_{2+iT}^{2+i\infty}+\int_{2-iT}^{2+iT}+
\int_{2-i\infty}^{2-iT})u^{-s}G(u;s)\frac{x^{s+1}}{s(s+1)}ds.
 $\\
We have \\
$ \sum\limits_{u\leq Q}\frac{1}{2\pi i}(\int_{2+iT}^{2+i\infty}+
\int_{2-i\infty}^{2-iT})u^{-s}G(u;s)\frac{x^{s+1}}{s(s+1)}ds\ll
\sum\limits_{u\leq Q}u^{-2}x^{3}T^{-1}\ll  x^{3}T^{-1},$\\
therefore
\bi
 S_{E}(x,Q)=\sum\limits_{u\leq Q}\frac{1}{2\pi
i}\int_{2-iT}^{2+iT}u^{-s}G(u;s)\frac{x^{s+1}}{s(s+1)}ds+O(x^{3}T^{-1})=A+O(x^{3}T^{-1}),
\ei
 where\\
$A=\sum\limits_{u\leq Q}\frac{1}{2\pi
i}\int_{2-iT}^{2+iT}u^{-s}G(u;s)\frac{x^{s+1}}{s(s+1)}ds$.\\
Let $\frac{1}{k}-1+\tau\leq\sigma<0$, by residue theorem \\
$ A=\sum\limits_{u\leq Q}\frac{1}{2\pi
i}\int_{\sigma-iT}^{\sigma+iT}u^{-s}G(u;s)\frac{x^{s+1}}{s(s+1)}ds
+\sum\limits_{u\leq Q}\frac{1}{2\pi
i}\int_{2-iT}^{\sigma-iT}u^{-s}G(u;s)\frac{x^{s+1}}{s(s+1)}ds\\
+\sum\limits_{u\leq Q}\frac{1}{2\pi
i}\int_{\sigma+iT}^{2+iT}u^{-s}G(u;s)\frac{x^{s+1}}{s(s+1)}ds
+{Res}_{s=0,1}(\sum\limits_{u\leq
Q}u^{-s}G(u;s)\frac{x^{s+1}}{s(s+1)}). $\\
{\bf LEMMA 10.1}\ For fixed $\frac{1}{k}-1+\tau\leq\sigma\leq 2,
\mid t\mid\geq 1$ and $s=\sigma+it$, we have the Stirling formula
$$\Gamma(s)=\sqrt{2\pi}e^{-\frac{\pi}{2}\mid t\mid}{\mid t\mid}^{\sigma-\frac{1}{2}}e^{it(\log \mid t\mid-1)}
e^{\lambda\frac{i\pi}{2}(\sigma-\frac{1}{2})}(1+O({\mid
t\mid}^{-1})),$$
where $\lambda=1 $ if  $t>0$, and $\lambda=-1 $ if
$t<0$; and
$$G(u;s)\ll \mid t\mid+\mid t\mid^{1-\sigma}, \ \ \ \ \frac{1}{k}-1+\tau\leq\Re s\leq 2, \ \mid t\mid\geq 1.$$
Proof.\ The first part is the Corollary 3.3.3 of [5], or see 4.42 of [6].\\
By Theorem 7.1.1 of [5], for $-1<\Re s\leq 2,\ y\geq 0$\\
$\zeta(s,a)=\sum\limits_{0\leq n\leq
y}\frac{1}{(n+a)^{s}}+\frac{(y+a)^{1-s}}{s-1}+\frac{y-[y]-\frac{1}{2}}{(y+a)^{s}}
-s\int_{y}^{\infty}
\frac{u-[u]-\frac{1}{2}}{(u+a)^{s+1}}du.$\\
When $\Re s\geq\frac{1}{2}$, we choose $y=1$, then\\
$\zeta(s,a)=a^{-s}+(1+a)^{-s}+\frac{(1+a)^{1-s}}{s-1}-\frac{1}{2}(1+a)^{-s}-s\int_{1}^{\infty}
\frac{u-[u]-\frac{1}{2}}{(u+a)^{s+1}}du,$\\
hence
$$\zeta(s,a)\ll a^{-\sigma}+\mid t\mid,\ \ \mid t\mid\geq 1 .$$
When $\frac{1}{k}-1<\Re s\leq \frac{1}{2}$, we choose $y=\mid t\mid\geq 1$, then\\
$\zeta(s,a)\ll a^{-\sigma}+\sum\limits_{1\leq n\leq \mid
t\mid}\frac{1}{(n+a)^{\sigma}}+\mid t\mid ^{1-\sigma}+\mid
s\int_{\mid t\mid}^{\infty}
\frac{u-[u]-\frac{1}{2}}{(u+a)^{s+1}}du\mid,$\\
by partial integrations,\\
$\mid s\int_{\mid t\mid}^{\infty}
\frac{u-[u]-\frac{1}{2}}{(u+a)^{s+1}}du\mid\ll \mid t\mid^{1-\sigma}+\mid t\mid^{2}\mid t\mid^{-\sigma-1}\ll \mid t\mid^{1-\sigma}.$\\
We have\\
$\zeta(s,a)\ll a^{-\sigma}+\mid t\mid+\mid t\mid^{1-\sigma}, \ \ \
\frac{1}{k}-1+\tau\leq\Re s\leq 2, \ \mid t\mid\geq 1.$\\
 By (9.6), (9.21) and $1\leq
A_{1},A_{3},A_{4},A_{5}\ll
d_{2}^{k}d_{3}^{k}$,\\
$$G(u;s)\ll 1+\mid t\mid+\mid t\mid^{1-\sigma}, \ \ \ \ \frac{1}{k}-1+\tau\leq\Re s\leq 2, \ \mid t\mid\geq 1,$$
the lemma follows.\\
By Lemma 10.1,\\
$\sum\limits_{u\leq Q}\frac{1}{2\pi
i}\int_{2-iT}^{\sigma-iT}u^{-s}G(u;s)\frac{x^{s+1}}{s(s+1)}ds ,\
\sum\limits_{u\leq Q}\frac{1}{2\pi
i}\int_{\sigma+iT}^{2+iT}u^{-s}G(u;s)\frac{x^{s+1}}{s(s+1)}ds\ll
x^{3}Q^{1-\sigma}T^{-\frac{1}{k}},$\\
therefore
\bi \left.
\begin{array}{cc}
 A= \sum\limits_{u\leq Q}\frac{1}{2\pi
i}\int_{\sigma-iT}^{\sigma+iT}u^{-s}G(u;s)\frac{x^{s+1}}{s(s+1)}ds\\
+{Res}_{s=0,1}(\sum\limits_{u\leq
Q}u^{-s}G(u;s)\frac{x^{s+1}}{s(s+1)})+O(x^{3}Q^{1-\sigma}T^{-\frac{1}{k}}).
\end{array}
\right. \ei
 Write \\
$ B_{1}=\sum\limits_{u\leq Q}\frac{1}{2\pi
i}\int_{\sigma-iT}^{\sigma+iT}u^{-s}G_{(1)}(u,s)\frac{x^{s+1}}{s(s+1)}ds,
$\\
then, by (9.19),\\
$ B_{1}=\sum\limits_{u\leq Q}\frac{1}{2\pi
i}\int_{\sigma-iT}^{\sigma+iT}u^{-s}\sum\limits_{\stackrel{(d_{2}d_{3},u)=1}{(d_{1}d_{2}d_{3},2)=1}}
h(d_{1},d_{2},d_{3})\frac{(u,2^{k}d_{1}^{k})^{s}}{2^{ks}d_{1}^{ks}}
\zeta(s;A_{1},d_{2}^{k}d_{3}^{k})\frac{x^{s+1}}{s(s+1)}ds\\
=\sum\limits_{u\leq Q}\frac{1}{2\pi
i}\int_{\sigma-iT}^{\sigma+iT}u^{-s}\sum\limits_{\stackrel{(d_{2}d_{3},u)=1}{(d_{1}d_{2}d_{3},2)=1}}
h(d_{1},d_{2},d_{3})\frac{(u,2^{k}d_{1}^{k})^{s}}{2^{ks}d_{1}^{ks}}
G(s)(d_{2}d_{3})^{-ks}(e^{\frac{1}{2}i\pi
s}\zeta(1-s;\xi_{d_{2}^{k}d_{3}^{k}}^{A_{1}})-e^{\frac{-1}{2}i\pi
s}\zeta(1-s;\xi_{d_{2}^{k}d_{3}^{k}}^{-A_{1}}))\frac{x^{s+1}}{s(s+1)}ds=B_{11}+B_{12},
$ say.\\

We have\\
 $
\int_{\sigma-iT}^{\sigma+iT}u^{-s}\frac{(u,2^{k}d_{1}^{k})^{s}}{2^{ks}d_{1}^{ks}}
G(s)(d_{2}d_{3})^{-ks}e^{\frac{1}{2}i\pi
s}\zeta(1-s;\xi_{d_{2}^{k}d_{3}^{k}}^{A_{1}})\frac{x^{s+1}}{s(s+1)}ds\\
=\sum\limits_{m=1}^{\infty}\int_{\sigma-iT}^{\sigma+iT}u^{-s}\frac{(u,2^{k}d_{1}^{k})^{s}}{2^{ks}d_{1}^{ks}}
G(s)(d_{2}d_{3})^{-ks}e^{\frac{1}{2}i\pi
s}\xi_{d_{2}^{k}d_{3}^{k}}^{A_{1}m}m^{-(1-s)}\frac{x^{s+1}}{s(s+1)}ds\\
=e^{\frac{1}{2}i\pi
\sigma}x^{\sigma+1}iu^{-\sigma}\frac{(u,2^{k}d_{1}^{k})^{\sigma}}{2^{k\sigma}d_{1}^{k\sigma}}
(d_{2}d_{3})^{-k\sigma}
\sum\limits_{m=1}^{\infty}\xi_{d_{2}^{k}d_{3}^{k}}^{A_{1}m}
m^{-1+\sigma}\int_{-T}^{T}u^{-it}\frac{(u,2^{k}d_{1}^{k})^{it}}{2^{kti}d_{1}^{kti}}
G(\sigma+it)(d_{2}d_{3})^{-kti}e^{-\frac{1}{2}\pi t}m^{it}\times\\
\frac{x^{it}}{(\sigma+it)(1+\sigma+it)}dt,$\\
and\\
$\int_{-T}^{T}u^{-it}\frac{(u,2^{k}d_{1}^{k})^{it}}{2^{kti}d_{1}^{kti}}
G(\sigma+it)(d_{2}d_{3})^{-kti}e^{-\frac{1}{2}\pi t}m^{it}
\frac{x^{it}}{(\sigma+it)(1+\sigma+it)}dt\\
\ll
(\int_{1}^{T}+\int_{-T}^{-1})u^{-it}\frac{(u,2^{k}d_{1}^{k})^{it}}{2^{kti}d_{1}^{kti}}
G(\sigma+it)(d_{2}d_{3})^{-kti}e^{-\frac{1}{2}\pi t}m^{it}
\frac{x^{it}}{(\sigma+it)(1+\sigma+it)}dt+1$.\\
Let $1\leq T_{1}\leq T_{2}\leq 2T_{1}$. By Lemma 10.1,\\
$\int_{T_{1}}^{T_{2}}u^{-it}\frac{(u,2^{k}d_{1}^{k})^{it}}{2^{kti}d_{1}^{kti}}
G(\sigma+it)(d_{2}d_{3})^{-kti}e^{-\frac{1}{2}\pi t}m^{it}
\frac{x^{it}}{(\sigma+it)(1+\sigma+it)}dt\\
=-i\int_{T_{1}}^{T_{2}}u^{-it}\frac{(u,2^{k}d_{1}^{k})^{it}}{2^{kti}d_{1}^{kti}}
(2\pi)^{\sigma+it-1}
\sqrt{2\pi}e^{-\frac{\pi}{2}t}t^{1-\sigma-\frac{1}{2}}e^{-it(\log
 t-1)}
e^{-\frac{i\pi}{2}(1-\sigma-\frac{1}{2})}(1+O(t^{-1}))\times\\
(d_{2}d_{3})^{-kti}e^{-\frac{1}{2}\pi t}m^{it}
\frac{x^{it}}{(\sigma+it)(1+\sigma+it)}dt\ll T_{1}^{-\frac{3}{2}-\sigma}, $\\
and\\
$\int_{T_{1}}^{T_{2}}u^{it}\frac{(u,2^{k}d_{1}^{k})^{-it}}{2^{-kti}d_{1}^{-kti}}
G(\sigma-it)(d_{2}d_{3})^{kti}e^{\frac{1}{2}\pi t}m^{-it}
\frac{x^{-it}}{(\sigma-it)(1+\sigma-it)}dt\\
=-i(2\pi)^{\sigma-\frac{1}{2}}e^{\frac{i\pi}{2}(\frac{1}{2}-\sigma)}\int_{T_{1}}^{T_{2}}u^{it}\frac{(u,2^{k}d_{1}^{k})^{-it}}{2^{-kti}d_{1}^{-kti}}
(2\pi)^{-it}t^{1-\sigma-\frac{1}{2}}e^{it(\log
 t-1)}(1+O(t^{-1}))\times\\
(d_{2}d_{3})^{kti}m^{-it}
\frac{x^{-it}}{(\sigma-it)(1+\sigma-it)}dt\\
=i(2\pi)^{\sigma-\frac{1}{2}}e^{\frac{i\pi}{2}(\frac{1}{2}-\sigma)}\int_{T_{1}}^{T_{2}}u^{it}\frac{(u,2^{k}d_{1}^{k})^{-it}}{2^{-kti}d_{1}^{-kti}}
(2\pi)^{-it}t^{-\frac{3}{2}-\sigma}e^{it(\log
 t-1)}
(d_{2}d_{3})^{kti}m^{-it}
x^{-it}dt+O(T_{1}^{-\frac{3}{2}-\sigma})$.\\
Write\\
$B=\int_{T_{1}}^{T_{2}}u^{it}\frac{(u,2^{k}d_{1}^{k})^{-it}}{2^{-kti}d_{1}^{-kti}}
(2\pi)^{-it}t^{-\frac{3}{2}-\sigma}e^{it(\log
 t-1)}
(d_{2}d_{3})^{kti}m^{-it}
x^{-it}dt=\int_{T_{1}}^{T_{2}}t^{-\frac{3}{2}-\sigma}e^{iF(t)}dt,$\\
where\\
$F(t)=t\log u-t\log(u,2^{k}d_{1}^{k})+kt\log
(2d_{1}d_{2}d_{3})-t\log (2\pi)+t(\log t-1)-t\log(mx).$\\

 $F^{\prime\prime}(t)=t^{-1}$,\\
therefore, by Lemma 4.4 of [7] and by partial integrations, \\
$\int_{T_{1}}^{T_{2}}e^{iF(t)}dt\ll T_{1}^{\frac{1}{2}}, \ \ B\ll
T_{1}^{-1-\sigma}$.\\
Hence
$$\int_{-T}^{T}u^{-it}\frac{(u,2^{k}d_{1}^{k})^{it}}{2^{kti}d_{1}^{kti}}
G(\sigma+it)(d_{2}d_{3})^{-kti}e^{-\frac{1}{2}\pi t}m^{it}\times\\
\frac{x^{it}}{(\sigma+it)(1+\sigma+it)}dt\ll 1,  $$
and, by (9.6),
$$B_{11}\ll \sum\limits_{u\leq Q}x^{\sigma
+1}u^{-\sigma}\ll x^{\sigma +1}Q^{1-\sigma}.$$
The same proof works for $B_{12}$, then \\
$$B_{12}\ll x^{\sigma +1}Q^{1-\sigma},$$
we obtain
$$\sum\limits_{u\leq Q}\frac{1}{2\pi
i}\int_{\sigma-iT}^{\sigma+iT}u^{-s}G_{(1)}(u,s)\frac{x^{s+1}}{s(s+1)}ds\ll
x^{\sigma +1}Q^{1-\sigma}.$$
Similarly, we deduce
$$\sum\limits_{u\leq Q}\frac{1}{2\pi
i}\int_{\sigma-iT}^{\sigma+iT}u^{-s}G_{(j)}(u,s)\frac{x^{s+1}}{s(s+1)}ds\ll
x^{\sigma +1}Q^{1-\sigma},\ \ j=3,4,5.$$
By (9.21),
 \bi
\sum\limits_{u\leq Q}\frac{1}{2\pi
i}\int_{\sigma-iT}^{\sigma+iT}u^{-s}G(u;s)\frac{x^{s+1}}{s(s+1)}ds\ll
x^{\sigma +1}Q^{1-\sigma}. \ei Let $T\rightarrow\infty$, by (10.1),
(10.2) and (10.3),
 \bi S_{E}(x,Q)={Res}_{s=0,1}(\sum\limits_{u\leq
Q}u^{-s}G(u;s)\frac{x^{s+1}}{s(s+1)})+O(x^{\sigma +1}Q^{1-\sigma}).
\ei We have
 \bi {Res}_{s=1}(\sum\limits_{u\leq
Q}u^{-s}G(u;s)\frac{x^{s+1}}{s(s+1)})=c_{2}(Q)x^{2}, \ei
 where $c_{2}(Q)=\frac{1}{2}\sum\limits_{u\leq
Q}u^{-1}{Res}_{s=1}G(u;s)$,  $c_{2}(Q)$ is independence of $x$.\\
We need the value $\zeta(0;a,k)=\frac{1}{2}-\frac{a}{k},\ 1\leq
a\leq  k$.\\
\bi {Res}_{s=0}(\sum\limits_{u\leq
Q}u^{-s}G(u;s)\frac{x^{s+1}}{s(s+1)})= x\sum\limits_{u\leq Q}G(u;0).
\ei
 {\bf LEMMA 10.2}\ We have\\
$\sum\limits_{u\leq
Q}G_{(1)}(u,0)=c_{31}Q+O(Q^{\frac{1}{k}+\varepsilon}),\ \
 \sum\limits_{u\leq Q}G_{(3)}(u,0)=c_{33}Q+O(Q^{\frac{1}{k}+\varepsilon}),\\
 \sum\limits_{u\leq
 Q}G_{(4)}(u,0)=c_{34}Q+O(Q^{\frac{1}{k}+\varepsilon}),\ \
 \sum\limits_{u\leq Q}G_{(5)}(u,0)=c_{35}Q+O(Q^{\frac{1}{k}+\varepsilon})$,\\
where $c_{31}, c_{33}, c_{34}, c_{35}$ are independence of $Q$.\\
Proof.\ We only consider the case $(d_{i}, d_{l})=1, i\not=l$. \\
\bi \left.
\begin{array}{cc}
 \sum\limits_{u\leq Q}G_{(1)}(u,0)
 = \sum\limits_{u\leq Q}\sum\limits_{\stackrel{(d_{2}d_{3},u)=1}{(d_{1}d_{2}d_{3},2)=1}}
h(d_{1},d_{2},d_{3})\zeta(0;A_{1},d_{2}^{k}d_{3}^{k})\\
=\sum\limits_{u\leq
Q}\sum\limits_{\stackrel{(d_{2}d_{3},u)=1}{(d_{1}d_{2}d_{3},2)=1}}
h(d_{1},d_{2},d_{3})(\frac{1}{2}-A_{1}d_{2}^{-k}d_{3}^{-k}).
\end{array}
\right. \ei
 We have \\
$\sum\limits_{\stackrel{u\leq Q}{(m,u)=1}}1
 =Q\sum\limits_{t\mid m}\mu(t)t^{-1}+O(\tau(m)),$\\
then\\
$\sum\limits_{u\leq
Q}\sum\limits_{\stackrel{(d_{2}d_{3},u)=1}{(d_{1}d_{2}d_{3},2)=1}}
h(d_{1},d_{2},d_{3})\frac{1}{2}
=\frac{1}{2}\sum\limits_{(d_{1}d_{2}d_{3},2)=1}h(d_{1},d_{2},d_{3})\sum\limits_{\stackrel{u\leq
Q}{(d_{2}d_{3},u)=1}}1\\
=\frac{1}{2}\sum\limits_{(d_{1}d_{2}d_{3},2)=1}h(d_{1},d_{2},d_{3})(Q\sum\limits_{t\mid
d_{2}d_{3}}\mu(t)t^{-1}+O(\tau(d_{2}d_{3}))),$\\
therefore, by (9.6),
 \bi
\sum\limits_{u\leq
Q}\sum\limits_{\stackrel{(d_{2}d_{3},u)=1}{(d_{1}d_{2}d_{3},2)=1}}
h(d_{1},d_{2},d_{3})\frac{1}{2}=\frac{1}{2}Q\sum\limits_{(d_{1}d_{2}d_{3},2)=1}h(d_{1},d_{2},d_{3})\sum\limits_{t\mid
d_{2}d_{3}}\mu(t)t^{-1}+O(1).
 \ei
 Write\\
$C=\sum\limits_{\stackrel{u\leq Q}{(d_{2}d_{3},u)=1}}A_{1}
=\sum\limits_{\stackrel{u\leq
Q}{(d_{2}d_{3},u)=1}}\{-b_{1}d_{3}^{k}\overline{d_{3}^{k}}+b_{2}d_{2}^{k}\overline{d_{2}^{k}}\}^{\prime},$\\
where
$A_{1}=\{-b_{1}d_{3}^{k}\overline{d_{3}^{k}}+b_{2}d_{2}^{k}\overline{d_{2}^{k}}\}^{\prime}$
means
$ 1\leq A_{1}\leq d_{2}^{k}d_{3}^{k}, (A_{1},d_{2}d_{3})=1.$\\
Let $(u,2^{k}d_{1}^{k})=j, \ u=ju_{1}$, then, by (9.8),(9.9) and (9.10),\\
$$ 2^{k}d_{1}^{k}u_{1}b_{1}\equiv 1 (mod\ d_{2}^{k}),\
   2^{k}d_{1}^{k}u_{1}b_{2}\equiv 1 (mod\ d_{3}^{k}), $$
and\\
$C=\sum\limits_{j_{1}=1}^{d_{2}^{k}}{}^{\prime}
\sum\limits_{j_{2}=1}^{d_{3}^{k}}{}^{\prime}
\{-j_{1}d_{3}^{k}\overline{d_{3}^{k}}+j_{2}d_{2}^{k}\overline{d_{2}^{k}}\}^{\prime}\sum\limits_{j\mid 2^{k}d_{1}^{k}}D$,\\
where\\
\[ \left.
\begin{array}{cc}
D=\#\{u_{1}: u_{1}\leq Qj^{-1}, (u_{1},2^{k}d_{1}^{k}j^{-1})=1,
2^{k}d_{1}^{k}u_{1}j_{1}\equiv
1 (mod\ d_{2}^{k}),\\
   2^{k}d_{1}^{k}u_{1}j_{2}\equiv 1 (mod\ d_{3}^{k})\}.
\end{array}
\right.
\]
$D=\sum\limits_{t\mid
2^{k}d_{1}^{k}j^{-1}}\mu(t)\sum\limits_{\stackrel{u\leq
Qj^{-1}t^{-1}}{ u2^{k}d_{1}^{k}tj_{1}\equiv 1 (mod\ d_{2}^{k}),
   u2^{k}d_{1}^{k}tj_{2}\equiv 1 (mod\ d_{3}^{k})}}1
   =\sum\limits_{t\mid
2^{k}d_{1}^{k}j^{-1}}\mu(t)\sum\limits_{\stackrel{u\leq
Qj^{-1}t^{-1}}{ u\equiv J (mod\ d_{2}^{k}d_{3}^{k})}}1,$\\
where\\
\bi
J \equiv
\overline{2^{k}d_{1}^{k}t}\overline{j_{1}}d_{3}^{k}\overline{d_{3}^{k}}+
\overline{2^{k}d_{1}^{k}t}\overline{j_{2}}d_{2}^{k}\overline{d_{2}^{k}}
 \ (mod\
d_{2}^{k}d_{3}^{k}),{\rm and } \ 1\leq J\leq d_{2}^{k}d_{3}^{k}. \ei
 We have
\[
\sum\limits_{\stackrel{u\leq Qj^{-1}t^{-1}}{ u\equiv J (mod\
d_{2}^{k}d_{3}^{k})}}1=\left\{\begin{array}{cc}
Qj^{-1}t^{-1}d_{2}^{-k}d_{3}^{-k}+O((Qj^{-1}t^{-1}d_{2}^{-k}d_{3}^{-k})^{\frac{1}{k}+\varepsilon}),&{\rm if}\ \ d_{2}^{k}d_{3}^{k}\leq Qj^{-1}t^{-1},\\
1,&{\rm if}\ \ d_{2}^{k}d_{3}^{k}> Qj^{-1}t^{-1}\geq J,\\
0,&{\rm if}\ \ d_{2}^{k}d_{3}^{k}\geq J> Qj^{-1}t^{-1} .
\end{array}
\right.
\]
 Therefore\\
$ \sum\limits_{u\leq
Q}\sum\limits_{\stackrel{(d_{2}d_{3},u)=1}{(d_{1}d_{2}d_{3},2)=1}}
h(d_{1},d_{2},d_{3})A_{1}d_{2}^{-k}d_{3}^{-k}\\
=\sum\limits_{(d_{1}d_{2}d_{3},2)=1}h(d_{1},d_{2},d_{3})d_{2}^{-k}d_{3}^{-k}
\sum\limits_{j_{1}=1}^{d_{2}^{k}}{}^{\prime}
\sum\limits_{j_{2}=1}^{d_{3}^{k}}{}^{\prime}
\{-j_{1}d_{3}^{k}\overline{d_{3}^{k}}+j_{2}d_{2}^{k}\overline{d_{2}^{k}}\}^{\prime}\sum\limits_{j\mid
2^{k}d_{1}^{k}}\sum\limits_{t\mid
2^{k}d_{1}^{k}j^{-1}}\mu(t)\sum\limits_{\stackrel{u\leq
Qj^{-1}t^{-1}}{ u\equiv J (mod\ d_{2}^{k}d_{3}^{k})}}1 \\
= \sum\limits_{(d_{1},2)=1}\sum\limits_{j\mid
2^{k}d_{1}^{k}}\sum\limits_{t\mid 2^{k}d_{1}^{k}j^{-1}}\mu(t)
\sum\limits_{d_{2}^{k}d_{3}^{k}\leq Qj^{-1}t^{-1}}
h(d_{1},d_{2},d_{3})d_{2}^{-k}d_{3}^{-k}\times\\
\sum\limits_{j_{1}=1}^{d_{2}^{k}}{}^{\prime}
\sum\limits_{j_{2}=1}^{d_{3}^{k}}{}^{\prime}
\{-j_{1}d_{3}^{k}\overline{d_{3}^{k}}+j_{2}d_{2}^{k}\overline{d_{2}^{k}}\}^{\prime}
(Qj^{-1}t^{-1}d_{2}^{-k}d_{3}^{-k}+O((Qj^{-1}t^{-1}d_{2}^{-k}d_{3}^{-k})^{\frac{1}{k}+\varepsilon}))\\
+ \sum\limits_{(d_{1},2)=1}\sum\limits_{j\mid
2^{k}d_{1}^{k}}\sum\limits_{t\mid 2^{k}d_{1}^{k}j^{-1}}\mu(t)
\sum\limits_{d_{2}^{k}d_{3}^{k}> Qj^{-1}t^{-1}}
h(d_{1},d_{2},d_{3})d_{2}^{-k}d_{3}^{-k}
\sum\limits_{j_{1}=1}^{d_{2}^{k}}{}^{\prime}
\sum\limits_{\stackrel{j_{2}=1}{1\leq J\leq
Qj^{-1}t^{-1}}}^{d_{3}^{k}}{}^{\prime}
\{-j_{1}d_{3}^{k}\overline{d_{3}^{k}}+j_{2}d_{2}^{k}\overline{d_{2}^{k}}\}^{\prime}.$\\
By (9.6), \\
$\sum\limits_{(d_{1},2)=1}\sum\limits_{j\mid
2^{k}d_{1}^{k}}\sum\limits_{t\mid 2^{k}d_{1}^{k}j^{-1}}\mu(t)
\sum\limits_{\stackrel{(d_{2}d_{3},2)=1}{d_{2}^{k}d_{3}^{k}\leq
Qj^{-1}t^{-1}}}
h(d_{1},d_{2},d_{3})d_{2}^{-k}d_{3}^{-k}\times\\
\sum\limits_{j_{1}=1}^{d_{2}^{k}}{}^{\prime}
\sum\limits_{j_{2}=1}^{d_{3}^{k}}{}^{\prime}
\{-j_{1}d_{3}^{k}\overline{d_{3}^{k}}+j_{2}d_{2}^{k}\overline{d_{2}^{k}}\}^{\prime}
(Qj^{-1}t^{-1}d_{2}^{-k}d_{3}^{-k}+O((Qj^{-1}t^{-1}d_{2}^{-k}d_{3}^{-k})^{\frac{1}{k}+\varepsilon}))\\
=Q\sum\limits_{(d_{1}d_{2}d_{3},2)=1}\sum\limits_{j\mid
2^{k}d_{1}^{k}}\sum\limits_{t\mid
2^{k}d_{1}^{k}j^{-1}}\mu(t)j^{-1}t^{-1}
h(d_{1},d_{2},d_{3})d_{2}^{-2k}d_{3}^{-2k}\sum\limits_{a=1}^{d_{2}^{k}d_{3}^{k}}{}^{\prime}a
+O(Q^{\frac{1}{k}+\varepsilon}),$\\
we note that : if $tj>Q$, this will produce the error term $O(1)$,
 therefore the above formula is valid.\\
By (9.6) and (10.9),\\
$\sum\limits_{(d_{1},2)=1}\sum\limits_{j\mid
2^{k}d_{1}^{k}}\sum\limits_{t\mid 2^{k}d_{1}^{k}j^{-1}}\mu(t)
\sum\limits_{d_{2}^{k}d_{3}^{k}> Qj^{-1}t^{-1}}
h(d_{1},d_{2},d_{3})d_{2}^{-k}d_{3}^{-k}\times\\
\sum\limits_{j_{1}=1}^{d_{2}^{k}}{}^{\prime}
\sum\limits_{\stackrel{j_{2}=1}{1\leq J\leq
Qj^{-1}t^{-1}}}^{d_{3}^{k}}{}^{\prime}
\{-j_{1}d_{3}^{k}\overline{d_{3}^{k}}+j_{2}d_{2}^{k}\overline{d_{2}^{k}}\}^{\prime}\\
\ll
\sum\limits_{(d_{1},2)=1}\sum\limits_{j\mid
2^{k}d_{1}^{k}}\sum\limits_{t\mid 2^{k}d_{1}^{k}j^{-1}}\mu(t)^{2}
\sum\limits_{d_{2}^{k}d_{3}^{k}> Qj^{-1}t^{-1}}
h(d_{1},d_{2},d_{3})d_{2}^{-k}d_{3}^{-k}d_{2}^{k}d_{3}^{k}Qj^{-1}t^{-1}\\
\ll Q\sum\limits_{(d_{1},2)=1}\sum\limits_{j\mid
2^{k}d_{1}^{k}}\sum\limits_{t\mid
2^{k}d_{1}^{k}j^{-1}}\mu(t)^{2}j^{-1}t^{-1}
\sum\limits_{d_{2}^{k}d_{3}^{k}> Qj^{-1}t^{-1}}
h(d_{1},d_{2},d_{3})\ll Q^{\frac{1}{k}+\varepsilon}.
 $\\
Hence, by (10.7) and (10.8),
$$\sum\limits_{u\leq
Q}G_{(1)}(u,0)=c_{31}Q+O(Q^{\frac{1}{k}+\varepsilon}), $$
where  $c_{31}$ is independence of $Q$.\\
Similarly, by the same method, we can obtain\\
$ \sum\limits_{u\leq
Q}G_{(3)}(u,0)=c_{33}Q+O(Q^{\frac{1}{k}+\varepsilon}),
 \sum\limits_{u\leq
 Q}G_{(4)}(u,0)=c_{34}Q+O(Q^{\frac{1}{k}+\varepsilon}), \\
 \sum\limits_{u\leq Q}G_{(5)}(u,0)=c_{35}Q+O(Q^{\frac{1}{k}+\varepsilon})$,\\
where $ c_{33}, c_{34}, c_{35}$ are independence of $Q$,  the lemma follows.\\
Consequently, for $\frac{1}{k}-1+\tau\leq\sigma<0$, by Lemma 10.2, (9.21), (10.4), (10.5) and (10.6),\\
$$ S_{E}(x,Q)=c_{2}(Q)x^{2}+c_{3}xQ+O(x^{\sigma
+1}Q^{1-\sigma}+xQ^{\frac{1}{k}+\varepsilon}), $$
by choosing
$\sigma=\frac{1}{k}-1+\tau <0$, we deduce
 \bi
S_{E}(x,Q)=c_{2}(Q)x^{2}+c_{3}xQ+O(x^{\frac{1}{k}+\tau}Q^{2-\frac{1}{k}-\tau}+xQ^{\frac{1}{k}+\varepsilon}),
\ei
 where $c_{2}(Q)$ is independence of $x$ , and $c_{3}$ is
independence of $x,Q$.
 \section{ Proof of the Theorem }
\setcounter{equation}{0}
 When $k=2$, the Theorem follows from Lemma 5.2. We suppose that $k>2$.\\
We also note that $R=x^{\frac{1}{2}+\tau}$. By (4.33), (5.6), (6.4), (7.2) and (10.9),\\
\bi S_{1}(x,Q)=2c_{2}(Q)x^{2}+2c_{3}xQ +[Q]\sum\limits_{n\leq
x}\mu_{k}(n)\mu_{k}(n+1)+O(x^{\frac{1}{k}+\tau}Q^{2-\frac{1}{k}-\tau}+
x^{\frac{3}{2}+\frac{1}{2k}+2\tau})\ei
 By Theorem 1 of [4],
$$
\sum\limits_{n\leq x}\mu_{k}(n)\mu_{k}(n+1)=\varrho
x+O(x^{\frac{2}{k+1}}),
$$
where $\varrho$ is defined by (4.17).\\
Hence, by (1.5), (1.6), (2.32), (3.4) and (11.1),\\
$$Y_{k}(x,Q)=c_{4}(Q)x^{2}+c_{5}xQ+O(x^{\frac{1}{k}+\tau}Q^{2-\frac{1}{k}-\tau}+
x^{\frac{3}{2}+\frac{1}{2k}+2\tau}),$$
 and
 \bi
c_{4}(Q)x^{2}+c_{5}xQ+O(x^{\frac{1}{k}+\tau}Q^{2-\frac{1}{k}-\tau}+
x^{\frac{3}{2}+\frac{1}{2k}+2\tau})
 \ll
 x^{\frac{2}{k}+\varepsilon}Q^{2-\frac{2}{k}}+x^{\frac{4}{k+1}+\varepsilon},
\ei
 where $c_{4}(Q)$ is independence of $x$ , and $c_{5}$ is
independence of $x,Q$. (11.2) is valid for all $1<Q\leq x$.\\
For fixed $Q$, we divide both sides of (11.2) by $x^{2}$ and let
$x\rightarrow \infty$, we deduce $c_{4}(Q)=0$ for all $Q>1$.\\
hence
$$Y_{k}(x,Q)=c_{5}xQ+O(x^{\frac{1}{k}+\tau}Q^{2-\frac{1}{k}-\tau}+
x^{\frac{3}{2}+\frac{1}{2k}+2\tau}),
 $$
and
 \bi
c_{5}xQ+O(x^{\frac{1}{k}+\tau}Q^{2-\frac{1}{k}-\tau}+
x^{\frac{3}{2}+\frac{1}{2k}+2\tau})
 \ll
 x^{\frac{2}{k}+\varepsilon}Q^{2-\frac{2}{k}}+x^{\frac{4}{k+1}+\varepsilon},
\ei We divide both sides of (11.3) by $xQ$ , by choosing
$Q=x^{\frac{5}{6}}$ and let $x\rightarrow \infty$, we obtain
$c_{5}=0$.\\
Consequently,
$$
Y_{k}(x,Q)=O( x^{\frac{1}{k}+\tau}Q^{2-\frac{1}{k}-\tau}+
x^{\frac{3}{2}+\frac{1}{2k}+2\tau}),
$$
since $\tau$ is a sufficiently small positive number, the Theorem
follows at once.\\

{\small School of Mathematical Sciences}\\
{\small Peking University}\\
{\small Beijing, 100871, P.R.China}   \\
{\small E-mail:\ mzzh@math.pku.edu.cn}
\end{document}